\documentclass[a4paper,11pt]{article}

\usepackage[T2A,T1]{fontenc}
\usepackage[russian,english]{babel}

\usepackage{a4wide}
\usepackage{amsmath}
\usepackage{amsfonts}
\usepackage{amssymb}
\usepackage{xypic}
\usepackage{url}

\newcommand\Reg{\mathbf{Reg}}
\newcommand\Sm[1][S]{\mathbf{Sm}/#1}
\newcommand\SmAff[1][S]{\mathbf{SmAff}_{#1}}
\newcommand\Ho[1][S]{\mathcal H(#1)}
\newcommand\Hopt[1][S]{\mathcal H_{\bullet}(#1)}
\newcommand\naive{\text{na\"ive}}
\newcommand\CH{\mathrm{ch}}
\newcommand\pf{\mathrm{pf}}
\newcommand\DM[1][S]{\mathrm{DM}(#1)}
\newcommand\SH[1][S]{\mathcal{SH}(#1)}
\newcommand\SHnaive[1][S]{\mathcal{SH}_{\naive}(#1)}
\newcommand\SHQ[1][S]{\mathcal{SH}(#1)_{\mathbf{Q}}}
\newcommand\SHpf[1][S]{\SH[#1]^{\pf}}
\newcommand\SHtop{\mathcal{SH}^{\mathrm{top}}}
\newcommand\SWpf[1][S]{SW(#1)^{\mathrm{ft}}}
\newcommand\Gm{\mathbf{G}_{\mathrm{m}}}
\newcommand\Sets{\mathbf{Sets}}
\newcommand\Ab{\mathbf{Ab}}
\newcommand\opp{\mathrm{opp}}
\newcommand\Nis{\mathrm{Nis}}
\newcommand\Presh[1][S]{{\Sm[#1]}^\opp\Sets}
\renewcommand\AA{\mathbf{A}}
\newcommand\NN{\mathbf{N}}
\newcommand\NNstar{\mathbf{N}-\{0\}}
\newcommand\ZZ{\mathbf{Z}}
\newcommand\Gr{\mathbf{Gr}}
\renewcommand\L{\mathbf{L}}
\newcommand\R{\mathbf{R}}
\newcommand\K{\mathbf{K}}
\newcommand\isomto{\overset{\sim}{\rightarrow}}
\newcommand\vers[1]{\overset{#1}{\rightarrow}}
\DeclareMathOperator{\rk}{rk}
\DeclareMathOperator{\Td}{Td}
\DeclareMathOperator{\Pic}{Pic}
\DeclareMathOperator{\End}{End}
\DeclareMathOperator{\Hom}{Hom}
\DeclareMathOperator{\Ext}{Ext}
\DeclareMathOperator{\Spec}{Spec}
\DeclareMathOperator{\SheafHom}{\mathbf{Hom}}
\DeclareMathOperator{\Sheafhom}{\mathbf{hom}}
\DeclareMathOperator\Th{Th}
\DeclareMathOperator{\oub}{oub}
\newcommand\eg{\emph{e.g.}}
\newcommand\ie{\emph{i.e.}}
\newcommand{\colim}{{\mathrm{colim}}}
\newcommand{\Rep}{\mathrm{R}}
\newcommand{\GL}{\mathbf{GL}}
\newcommand{\Univ}[1]{\mathrm{Univ}_{#1}}
\newcommand{\id}{\mathrm{id}}
\newcommand{\ID}[1]{\underline{\mathrm{id}}_{#1}}
\newcommand{\BGL}{\mathbf{BGL}}
\newcommand{\BGLS}[1][S]{\mathbf{BGL}_{#1}}
\newcommand{\BGLQ}{\mathbf{BGL}_{\mathbf{Q}}}
\newcommand{\BGLQS}[1][S]{\mathbf{BGL}_{\mathbf{Q},#1}}
\newcommand{\BGLnaive}{\mathbf{BGL}_{naive}}
\newcommand\SpecZ{\Spec(\mathbf{Z})}

\newcommand{\HBeilinson}{\mathbf{H}_{%
\text{\selectlanguage{russian}\emph{\CYRB}}}}
\newcommand{\HMot}[1][\mathbf{Z}]{\mathbf{H}_{#1}}
\newcommand{\HMotQ}{\HMot[\mathbf{Q}]}
\newcommand{\Htop}[1][\mathbf{Z}]{\mathbf{H}_{#1}^{\mathrm{top}}}

\newtheorem{theorem}[subsubsection]{Theorem}
\newtheorem{definition}[subsubsection]{Definition}
\newtheorem{lemma}[subsubsection]{Lemma}
\newtheorem{proposition}[subsubsection]{Proposition}
\newtheorem{remark}[subsubsection]{Remark}
\newtheorem{example}[subsubsection]{Example}
\newtheorem{corollary}[subsubsection]{Corollary}
\newtheorem{exercise}[subsubsection]{Exercise}

\newtheorem{theorem1}[subsection]{Theorem}

\newtheorem{theorem3}[paragraph]{Theorem}
\newtheorem{definition3}[paragraph]{Definition}
\newtheorem{lemma3}[paragraph]{Lemma}
\newtheorem{proposition3}[paragraph]{Proposition}
\newtheorem{remark3}[paragraph]{Remark}

\newtheorem{corollary3}[paragraph]{Corollary}

\title{Algebraic $K$-theory, $\mathbf{A}^1$-homotopy and Riemann-Roch theorems}
\author{Jo\"el Riou\;\footnote{Adresses: \texttt{joel.riou@math.u-psud.fr};
Universit\'e Paris-Sud 11, D\'epartement de
math\'ematiques, b\^at. 425, 91405 Orsay, France.}}
\date{September 9, 2009}

\begin{document}

\maketitle

\begin{abstract}
In this article, we show that the combination of the constructions done in
SGA~6~\cite{SGA6} and the $\mathbf{A}^1$-homotopy theory~\cite{MV}
naturally leads to
results on higher algebraic $K$-theory. This applies to the operations on
algebraic $K$-theory, Chern characters and Riemann-Roch theorems.
\end{abstract}

\tableofcontents

\bigskip


The starting point of this article is the theorem which represents the
algebraic $K$-theory of regular schemes in the $\mathbf{A}^1$-homotopy
theory:

\begin{theorem1}[{Morel-Voevodsky
\cite[Theorem~3.13, page 140]{MV}}]\label{theorem-morel-voevodsky}
\sloppy Let $S$ be a regular scheme.
Then, for any $n\in\mathbf{N}$ and $X\in\Sm$,
there is a canonical isomorphism
\[\Hom_{\Hopt}(S^n\wedge X_+,\mathbf{Z}\times\Gr)\simeq K_n(X)\;\text{.}\]
\end{theorem1}

Here, $\Gr$ is the colimit of the system
$(\Gr_{d,r})_{(d,r)\in\mathbf{N}^2}$ in the category of 
presheaves over $\Sm_\Nis$ where $\Gr_{d,r}$ is the Grassmann scheme
which parametrises subbundles of rank $d$ in the trivial bundle of rank
$d+r$. To make the definition of the transition morphisms unambiguous
enough, we may say that they are of the form $\Gr_{d,r}\to\Gr_{1+d,r}$ and
$\Gr_{d,r}\to\Gr_{d,r+1}$ and that the place where ``$1$'' appears tells us
on which side a trivial bundle of rank $1$ is added.

It should be pointed out that theorem~\ref{theorem-morel-voevodsky} only
applies to regular schemes because it can be true only over schemes where
the algebraic $K$-theory is known to be $\mathbf{A}^1$-invariant. This is
the reason why the assumption that the base scheme is regular will appear
throughout the paper.

From theorem~\ref{theorem-morel-voevodsky}, it follows that the
endomorphisms of $\mathbf{Z}\times\Gr$ in $\Hopt$ act on all the algebraic
$K$-groups of schemes in $\Sm$. The basic result we obtain in
section~\ref{section-unstable-results} is that these endomorphisms are
completely characterised by their action on $K_0$:

\begin{theorem1}\label{theorem-endomorphisms-of-z-gr}
Let $S$ be a regular scheme. We let $K_0(-)$ be the presheaf of sets on
$\Sm$ which maps $X$ to $K_0(X)$. Then, the map induced by
theorem~\ref{theorem-morel-voevodsky} is a bijection:
\[\End_{\Ho}(\mathbf{Z}\times\Gr)\isomto
\End_{\Sm^\opp\Sets}(K_0(-))\;\text{,}\]
where $\Sm^\opp\Sets$ is the category of presheaves of sets on $\Sm$.
\end{theorem1}

It follows that the operations defined in \cite{SGA6} at the level of $K_0$
(\eg, $\lambda^n$, $\Psi^k$) uniquely lift in $\Ho$. From there, using
theorem~\ref{theorem-morel-voevodsky}, we can make them act on higher
algebraic $K$-theory. This principle also works for operations involving
several operands (\eg, products) and in a sense which will be made precise
in section~\ref{section-algebraic-structures}, we obtain a machinery which
takes as an input the algebraic structures on $K_0$ and outputs such a
structure on $\mathbf{Z}\times\Gr$ inside $\Ho$. Thus,
$\mathbf{Z}\times\Gr$ is equipped with a structure of special
$\lambda$-ring with duality.

Structures of (special) $\lambda$-ring had already been obtained on higher
$K$-theory, with different scales of generality. We may mention
constructions of products, $\lambda$-operations or Adams operations
by Loday \cite{loday},
Waldhausen \cite{waldhausen}, Kratzer \cite{kratzer}, Soul\'e \cite{soule},
Grayson \cite{grayson}, Lecomte \cite{lecomte} and Levine \cite{levine-lambda}. We compare the
structures on $K_\star(X)$ for $X$ regular obtained by our method to these
previous constructions in section~\ref{section-comparison}. The comparison
with Waldhausen's product (see
proposition~\ref{proposition-comparison-waldhausen}) may seen surprisingly
straightforward, but it is a typical use of
theorem~\ref{theorem-endomorphisms-of-z-gr} and its variants involving
several operands (see
theorem~\ref{theorem-operations-k-theory-several-operands}).

Section~\ref{section-virtual-categories} relates our results to virtual
categories, an insight of Deligne \cite{deligne-determinant}. We show that,
after inverting $2$, constructions done at the level of $K_0$ refine to
these virtual categories, which embodies both $K_0$ and $K_1$. This theory
was used by Dennis Eriksson in his thesis \cite{eriksson-these}
in order to refine Riemann-Roch
theorems at the level of these virtual categories.

In section~\ref{section-additive-and-stable-results}, we focus on operations
$\tau\colon K_0(-)\to K_0(-)$ such that $\tau(x+y)=\tau(x)+\tau(y)$, \ie,
$H$-group endomorphisms of $\mathbf{Z}\times\Gr$ in $\Hopt$. We compute them
using the splitting principle. We show that the datum of $\tau$ is equivalent
to the datum of an element in $K_0(S)[[U]]$. Then, we construct, up to a unique
isomorphism in the stable homotopy category $\SH$, the $\mathbf{P}^1$-spectrum
$\BGL$ which represents algebraic $K$-theory and study its endomorphisms (it is
somewhat related but quite different from the methods of
\cite[Chapter~6]{adams}, \cite{adams-clarke} and \cite{anderson}). After
tensoring with $\mathbf{Q}$, we show that this spectrum decomposes in $\SH$ as
the direct sums of ``eigenspaces'' for the Adams operations. Alternate
interesting descriptions of stable operations on algebraic $K$-theory (and
more general oriented theory) have been obtained by very different methods by
Naumann, \O{}stv\ae{}r and Spitzweck in \cite{naumann-ostvaer-spitzweck}.

We prove in section~\ref{section-riemann-roch} that these ideas can be used
to obtain an homotopical variant of some Riemann-Roch theorems in the case
of a smooth and projective morphism $f\colon X\to S$. Basically, we prove
that certain Riemann-Roch formulas are satisfied on zeroth $K$-groups if
and only if they are satisfied on the whole higher algebraic $K$-theory.
In that section, we give formulas for the group of morphisms $\BGL\to
\HMot[A][n]$ in $\SH[k]$ where $k$ is a perfect field and $\HMot[A]$ the
motivic Eilenberg-Mac Lane spectrum with coefficients in $A$. This
computation gives a simple example of nonzero stably phantoms morphisms in
the $\mathbf{P}^1$-stable homotopy category $\SH[k]$: all morphisms
$\BGL\to\HMot{}[1]$ are stably phantoms. There is an homologous computation
in the standard topological stable homotopy category: this gives a more
concrete example than the one constructed in
\cite[Proposition~6.10]{christensen}.

If section~\ref{section-riemann-roch} stands as a significant exception,
most of these results appeared in my thesis \cite{riou-these} and were
announced in \cite{riou-cras-operations} (however, when different proofs
were available, my choices have tended to be different).
Hence, I would like to thank 
Yves Andr\'e,
Joseph Ayoub,
Denis-Charles Cisinski,
Fr\'ed\'eric D\'eglise,
Dennis Eriksson,
Hinda Hamraoui,
Bruno Kahn, 
Florence Lecomte,
Georges Maltsiniotis,
Fabien Morel,
Christophe Soul\'e,
Burt Totaro,
J\"org Wildeshaus
for their useful comments or discussions.

\section{First unstable results}
\label{section-unstable-results}
\subsection{Statements}

In this paper, we shall say that a scheme is regular if it is noetherian
separated and that all its local rings are regular local rings (see
\cite[IV~\S{}D]{ALM}). For any scheme $S$, the category of smooth and separated schemes of finite type over $S$ is denoted $\Sm$.

For regular schemes, all the standard definitions of algebraic $K$-theory
agree. Then, we may define some objects in the category $\Presh$ of presheaves
of sets over $\Sm$: for any natural number $n$, the presheaf that maps
$X\in\Sm$ to its $n$th algebraic $K$-group $K_n(X)$ is denoted $K_n(-)$.

\begin{theorem}\label{theorem-main}
Let $S$ be a regular scheme. For any natural
transformation $\tau\colon K_0(-)\to K_0(-)$ of presheaves of
sets on $\Sm$ such that $\tau(0)=0$, there is a unique reasonable way to
define an extension of $\tau\colon K_n(-)\to K_n(-)$ for all $n$.
\end{theorem}

This theorem is a consequence of the following $\AA^1$-homotopy theoretic
statement:

\begin{theorem}\label{theorem-endomorphisms-z-gr}
Let $S$ be a regular scheme. Then, the canonical map induced by the
isomorphism of theorem~\ref{theorem-morel-voevodsky} is a bijection:
\[\End_{\Ho}(\ZZ\times\Gr)\isomto \End_{\Presh}(K_0(-))\;\text{.}\]
\end{theorem}

Indeed, if $\tau\colon K_0(-)\to K_0(-)$ is a natural transformation, the
theorem says that there exists a unique morphism $\tilde{\tau}\colon
\ZZ\times\Gr\to\ZZ\times\Gr$ in $\Ho$ inducing $\tau$ on $K_0(-)$. As
$\ZZ\times\Gr$ has a structure of $H$-group (see \cite[page~139]{MV}), if
we assume $\tau(0)=0$, then we see that $\tilde{\tau}$ can be
identified to an endomorphism of $\ZZ\times\Gr$ in $\Hopt$. Such
endomorphisms not only induce natural transformations on $K_0(-)$ but also
on $K_n(-)$ for all $n$ as one may evaluate them on higher homotopy groups.

This theorem applies to operations like the $\lambda$-operations
$\lambda^n$ for all $n\in\NN$ \cite[V~2.2~b]{SGA6}, $\gamma$-operations
$\gamma^n$ for all $n\in\NNstar$ \cite[V~3.2]{SGA6} and Adams operations
$\Psi^k$ \cite[V~7.1]{SGA6} for all $k\in\mathbf{Z}$. Then, to construct
these operations on higher $K$-groups, the only specific information we
need to know is how to define them on $K_0$, which is usually easy using
the presentation of these groups by generators and relations. 

\begin{remark}
\label{remark-variant-gr}
One can prove similar results for $\Gr$ instead of $\mathbf{Z}\times\Gr$:
endomorphisms of $\Gr$ in $\Ho$ identify to endomorphisms of
$\tilde{K}_0(-)$ in $\Sm^\opp\Sets$ where $\tilde{K}_0(X)$ is the kernel of
the rank map $K_0(X)\to \mathbf{Z}^{\pi_0(X)}$. Moreover, in the situation
of theorem~\ref{theorem-main}, if we use the fact that the loop space $\R\Omega
(\mathbf{Z})$ of $\mathbf{Z}$ is $\bullet$,
we see that $\tau\colon K_n(-)\to K_n(-)$ for $n\geq
1$ only depends on the natural transformation $\tilde{K}_0(-)\to
\tilde{K}_0(-)$ induced by $\tau\colon K_0(-)\to K_0(-)$.
\end{remark}

The operations considered above are unary operations on algebraic
$K$-theory. One may also consider operations involving several operands
(\eg, the product law $K_0(X)\times K_0(X)\to K_0(X)$):

\begin{theorem}\label{theorem-operations-k-theory-several-operands}
Let $S$ be a regular scheme. Let $n$ be a natural number. Then, the
canonical map is a bijection:
\[\Hom_{\Ho}((\ZZ\times\Gr)^n,\ZZ\times\Gr)\to\Hom_{\Presh}(K_0(-)^n,K_0(-))\;\text{.}\]
\end{theorem}

As we shall see, the method of the proof allows to consider not only
operations on algebraic $K$-theory but also maps from algebraic $K$-theory
to other cohomology theories. However, we need to know that the cohomology
theory is represented by an object in $\Hopt$, which means that it can be
expressed as homotopy presheaves of an object in $\Hopt$:

\begin{definition}\label{definition-pi-zero}
Let $S$ be a noetherian scheme. Let $E$ be an object in $\Ho$. We let
$\pi_0E$ be the presheaf of sets on $\Sm$ defined by
$\pi_0E(X)=\Hom_{\Ho}(X,E)$. If $E$ belongs to $\Hopt$ and $n$ is any
natural number, we define a presheaf $\pi_nE$
by the formula $\pi_nE(X)=\pi_0\R\Omega^n E$, where $\R\Omega\colon
\Hopt\to \Hopt$ is the loop space functor.
\end{definition}

Theorem~\ref{theorem-morel-voevodsky} states that for any natural number
$n$ and $S$ a regular scheme, we have a canonical isomorphism
$\pi_n(\ZZ\times\Gr)\simeq K_n(-)$ in $\Presh$.

\begin{theorem}\label{theorem-property-k-implies-computation}
Let $S$ be a regular scheme. Let $E$ be an object in $\Hopt$. If we assume
that $E$ satisfies property~(K) (a mild technical assumption,
see definition~\ref{definition-property-k}), then the canonical map is a
bijection:
\[\Hom_{\Ho}(\ZZ\times\Gr,E)\isomto \Hom_{\Presh}(K_0(-),\pi_0E)\;\text{.}\]
This set of morphisms can also
be identified to an infinite product indexed by
$\mathbf{Z}$ of copies of the projective limit $\lim_{(d,r)\in\mathbf{N}^2}
\Hom_{\Ho}(\Gr_{d,r},E)$.
\end{theorem}

As above, there is a similar homotopical description of natural
transformations $K_0(-)^n\to \pi_0E$ involving $n$ operands.

We may focus on the $1$-operand case. If a natural transformation
$\tau\colon K_0(-)\to \pi_0E$ verifies $\tau(0)=0$, it corresponds to a
unique morphism $\ZZ\times\Gr\to E$ in $\Hopt$. Then, in the same way we
mentioned it for operations on algebraic $K$-theory, $\tau$ will induce
natural transformations $\tau\colon K_n(-)\to \pi_nE$ for all $n$.

The proof of theorems~\ref{theorem-endomorphisms-z-gr},
\ref{theorem-operations-k-theory-several-operands} and
\ref{theorem-property-k-implies-computation} will also supply a
concrete computation of the set of all operations on algebraic $K$-theory.
In the $1$-operand case, it gives:

\begin{theorem}\label{theorem-computation-k-of-z-gr}
Let $S$ be a regular scheme. The sets of endomorphisms
$\End_{\Ho}(\ZZ\times\Gr)\simeq \End_{\Presh}(K_0(-))$ can be identified to
the product $R^\ZZ$ of an infinite number of copies of a ring
$R=K_0(S)[[\tilde{\gamma}^1,\tilde{\gamma}^2,\dots]]$ of formal power
series with an
infinite number of variables and coefficient ring $K_0(S)$. The elements
$\tilde{\gamma}^n$ are related to the usual $\gamma$-operations on algebraic
$K$-theory.
\end{theorem}

The computation of the set of morphisms $\mathbf{Z}\times\Gr\to
\R\Omega^i(\mathbf{Z}\times\Gr)$ in $\Ho$ is given by a similar formula,
where $K_0(S)$ is replaced by $K_i(S)$.

\subsection{Proofs}

\begin{lemma}\label{lemma-milnor-exact-sequence}
Let $S$ be a noetherian scheme. Let $E$ be a group object in
$\Hopt$ (\ie, $E$ is an $H$-group). Let $(X_i)_{i\in\mathcal I}$ be a
direct system indexed by a directed ordered set $\mathcal I$. The colimit
of this system in the category of presheaves over $\Sm$
is denoted $\mathcal X$.
We assume that $\mathcal I$ has a cofinal sequence (\ie, there exists a
functor $x\colon \mathbf{N}\to \mathcal I$ such that for any $i\in I$,
there exists $n\in\mathbf{N}$ such that $i\leq x_n$).
Then, there is an exact sequence
of groups.
\[1\to\R^1\underset{i\in\mathcal I} {\lim}\; \pi_1E(X_i)\to
\Hom_{\Ho}(\mathcal X,E)\to \underset{i\in\mathcal I}\lim\; \pi_0 E(X_i)\to
1\;\text{.}\]
\end{lemma}

Using a cofinal sequence $\NN\to\mathcal I$, one may assume that $\mathcal
I=\NN$. In that case, it follows from the usual Milnor exact sequence
\cite[Proposition~VI.2.15]{SHT}.

\begin{definition}\label{definition-property-k}
With the notations of lemma~\ref{lemma-milnor-exact-sequence}, we say that
the direct system $(X_i)_{i\in\mathcal I}$ does not unveil phantoms in $E$
if the group $\R^1\underset{i\in\mathcal I}{\lim}\;\pi_1 E(X_i)$ vanishes.
We say that $E$ satisfies property (K) if the direct system
$(\Gr_{d,r})_{(d,r)\in\NN^2}$ does not unveil phantoms in $E$. More
generally, for any natural number $n$, we say that $E$ satisfies property
(K) with $n$ operands if the direct system $\left(\prod_{i=1}^n
\Gr_{d_i,r_i}\right)_{(d_1,r_1,\dots,d_n,r_n)\in\NN^{2n}}$ does not unveil
phantoms in $E$.
\end{definition}

Thus, whenever an inductive system $(X_i)_{i\in\mathcal I}$ does not unveil
phantoms in $E$, the datum of a morphism $\underset{i\in\mathcal
I}\colim\;X_i\to E$ in $\Ho$ is equivalent to the datum of a compatible
family of morphisms $X_i\to E$ in $\Ho$.

\begin{definition}
We let $\mathcal T$ be the family of morphisms in $\Sm$ of the form $T\to
X$ where $T$ is a torsor under a vector bundle over $X$.
\end{definition}

Locally on the base, morphisms in $\mathcal T$ are of the form $\AA^n\times
X\to X$. This implies that they induce $\AA^1$-weak equivalences. The
important fact we need about this family of maps is:

\begin{theorem}[{Jouanolou \cite[Lemme~1.5]{jouanolou}, Thomason
\cite[Proposition~4.4]{weibel-kh}}]\label{theorem-jouanolou-s-trick}~

Let $S$ be a regular scheme. For any $X\in\Sm$, there exists a morphism
$T\to X$ in $\mathcal T$ such that $T$ is an affine scheme.
\end{theorem}

We require that the scheme $T$ is affine; as $S$ is separated, it implies
that $T\to S$ in an affine morphism, but the converse implication is not
true. In the sequel, the word ``affine'' will be used in that absolute
sense only.

\begin{definition}\label{definition-generation-up-to-t}
Let $S$ be a regular scheme. Let $\mathcal X$ be a presheaf of sets on
$\Sm$. Then $\mathcal X$ defined an object in $\Ho$ and a presheaf of sets $\pi_0\mathcal X$ is attached to it.
We say that $\pi_0\mathcal X$ is generated by $\mathcal X$ up to
$\mathcal T$ if for any affine scheme $U\in\Sm$, the map $\mathcal
X(U)\to\pi_0\mathcal X(U)$ is onto.
\end{definition}

We will give an explanation for this terminological choice in
remark~\ref{remark-explanation-generation-up-to-t}. First, we see how one
may apply this definition to algebraic $K$-theory:

\begin{lemma}
\label{lemma-generation-up-to-t-of-z-gr}
Let $S$ be a regular scheme. If $\mathcal X=\ZZ\times\Gr$, then
$\pi_0\mathcal X$ is generated by $\mathcal X$ up to $\mathcal T$.
The same conclusion applies to $(\ZZ\times\Gr)^n$ for any
natural number $n$ and also to $(\mathbf{P}^\infty)^n$.
\end{lemma}

Obviously, the condition we have to check is stable under finite products.
Then, we shall first focus on
the case $\mathcal X=\ZZ\times\Gr$. It is implicit in the
proof of theorem~\ref{theorem-morel-voevodsky} that for any $n\in\ZZ$ and
$(d,r)\in\NN^2$, if we consider the canonical inclusion
$\iota_{d,r,n}
\colon\Gr_{d,r}=\{n\}\times\Gr_{d,r}\to \ZZ\times\Gr$ as an element in $\mathcal
X(\Gr_{d,r})$, its image in $\pi_0\mathcal X(\Gr_{d,r})$ corresponds to the
class $[\mathcal M'_{d,r}]-d+n$ in $K_0(\Gr_{d,r})$ under the isomorphism
of theorem~\ref{theorem-morel-voevodsky}, where $\mathcal M'_{d,r}$ is the
universal vector bundle of rank $d$ on $\Gr_{d,r}$. Then, the lemma follows
from the obvious fact that if $U$ is a connected affine scheme in $\Sm$,
any class $x\in K_0(U)$ is of the form $x=[\mathcal M]-d+n$ for some
integers $d$, $n$, and $\mathcal M$ a vector bundle of rank $d$ on $U$.
Indeed, as $U$ is affine, $\mathcal M$ is isomorphic to
a direct factor of $\mathcal
O_U^{d+r}$ for a big enough $r$. Then, by definition of Grassmann
varieties, there exists an $S$-morphism $f\colon U\to \Gr_{d,r}$ such that
$f^\star \mathcal M'_{d,r}\simeq \mathcal M$. It 
follows that the element in $\mathcal X(U)$
corresponding to the composition $\iota_{d,r,n}\circ f\colon U\to \mathcal
X$ maps to $x=f^\star([\mathcal M'_{d,r}-d+n])$ in
$\pi_0 \mathcal X(U)\simeq K_0(U)$.

The case $\mathcal X=\mathbf{P}^\infty$ is similar: it uses the identification $\pi_0\mathbf{P}^\infty=\Pic(-)$, see
\cite[Proposition~3.8, page 138]{MV}.

\begin{remark}\label{remark-explanation-generation-up-to-t}
The category of presheaves on $\Sm$ contains the full
subcategory of the category of
presheaves $\mathcal X$ such that for any $f\colon T\to X$
in $\mathcal T$, the map $f^\star\colon \mathcal X(X)\to\mathcal X(T)$ is a
bijection. This subcategory can be identified to the category of presheaves
on the localised category $\Sm{}[\mathcal T^{-1}]$ (see
\cite[Lemma~I.1.2]{gabriel-zisman}). For any presheaf
$\mathcal X$ on $\Sm$, there exists a universal presheaf $\mathcal
X[\mathcal T^{-1}]$ on
$\Sm{}[\mathcal T^{-1}]$ equipped with a morphism $\mathcal X\to
\mathcal X[\mathcal T^{-1}]$ (see \cite[I~5.1]{SGA4Vol1}).
As $\pi_0\mathcal X$ factors through
$\Sm{}[\mathcal T^{-1}]$, the canonical morphism $\mathcal X\to\pi_0\mathcal
X$ induces a morphism $\mathcal X[\mathcal T^{-1}]\to \pi_0\mathcal X$.
Using theorem~\ref{theorem-jouanolou-s-trick}, it
is easy to check that the condition stated in 
definition~\ref{definition-generation-up-to-t} implies that $\mathcal
X[\mathcal T^{-1}]\to \pi_0\mathcal X$ is an epimorphism. The converse
implication is also true, but we will not need it in the sequel.
This is the reason why we chose to refer to 
``generation up to $\mathcal T$'' in the terminology.

Moreover, the proof of lemma~\ref{lemma-generation-up-to-t-of-z-gr}
actually shows that as a presheaf $F$ on $\Sm{}[\mathcal T^{-1}]$
satisfying $F(X\sqcup Y)\isomto F(X)\times F(Y)$ for all $X$ and $Y$ in
$\Sm$, $K_0(-)\simeq
\pi_0(\mathbf{Z}\times\Gr)$ is generated by the elements $u_{d,r}+n$ for
all $(d,r)\in\mathbf{N}^2$ and $n\in\mathbf{Z}$.
\end{remark}

\begin{remark}
If is easy to deduce from theorem~\ref{theorem-jouanolou-s-trick} that the
localised category $\Sm{}[\mathcal T^{-1}]$ is equivalent to
$\SmAff{}[\mathcal H_{\mathbf{A}^1}^{-1}]$ where $\SmAff$ is the
fullsubcategory of $\Sm$ consisting of affine schemes and $\mathcal
H_{\mathbf{A}^1}$ is the family of projections $X\times \mathbf{A}^1\to X$
for $X\in\SmAff$ (see \cite[\S7.4]{kahn-sujatha}).
Hence, the category of $\mathcal T$-invariant presheaves
on $\Sm$ is equivalent to the category of $\mathbf{A}^1$-invariant
presheaves on $\SmAff$.
\end{remark}

\begin{proposition}\label{proposition-key}
Let $S$ be a regular scheme. Let $E\in\Hopt$ be an $H$-group. Let
$(X_i)_{i\in\mathcal I}$ be a direct system in $\Sm$
that does not unveil phantoms in $E$. We let $\mathcal X$ be the colimit of
this system in the category of presheaves over $\Sm$.
We assume that $\pi_0\mathcal X$ is generated by
$\mathcal X$ up to $\mathcal T$.
Then, the following obvious maps are bijections:
\[
\xymatrix{\Hom_{\Ho}(\mathcal X,E)\ar[r]^-\sim_-\alpha
\ar[rd]^{\sim}_\gamma&
\Hom_{\Presh}(\pi_0\mathcal X,\pi_0 E) \ar[d]^{\sim}_\beta\\
& \underset{i\in\mathcal I}{\lim}\;\pi_0 E(X_i)
}
\]
\end{proposition}

Using lemma~\ref{lemma-milnor-exact-sequence}, we see that the assumption
that $(X_i)_{i\in\mathcal I}$ does not unveil phantoms on $E$ precisely
says that $\gamma$ is a bijection. To finish the proof, we only have to
prove that $\beta$ is an injection. To do this, we may observe that 
$\underset{i\in\mathcal I}{\lim}\;\pi_0 E(X_i)$ identifies to
$\Hom_{\Presh}(\mathcal X,\pi_0 E)\simeq \Hom_{\Presh}(\mathcal X[\mathcal
T^{-1}],\pi_0E)$. Then, $\beta$ identifies to the map obtained by applying
the functor $\Hom_{\Presh}(-,\pi_0E)$ to the canonical map $\mathcal
X[\mathcal T^{-1}]\to \pi_0\mathcal X$, which is an epimorphism as
$\mathcal X\to \pi_0\mathcal X$ is an epimorphism up to $\mathcal T$ (see
remark~\ref{remark-explanation-generation-up-to-t}). Thus, $\beta$ is
injective.

\bigskip

At this stage, theorem~\ref{theorem-property-k-implies-computation} is
proved as lemma~\ref{lemma-generation-up-to-t-of-z-gr} implies that
it is a special case of proposition~\ref{proposition-key}. To
finish the proof of theorems \ref{theorem-main},
\ref{theorem-endomorphisms-z-gr} and
\ref{theorem-operations-k-theory-several-operands}, the remaining step is
the following lemma:

\begin{lemma}\label{lemma-property-k-for-z-gr}
Let $S$ be regular scheme. Let $n$ be a natural number.
The object $\ZZ\times\Gr$ satisfies property (K) with $n$ operands.
This conclusion also applies to the loop spaces $\R\Omega^j(\ZZ\times\Gr)$
for any $j\in\NN$.
\end{lemma}

On the one hand we have to notice the technical fact
that $\ZZ\times\Gr$ has a structure of $H$-group (see \cite[page 139]{MV}).
On the other hand,
we have to prove the vanishing of the $\R^1\lim$ of some projective
systems. To do this, one may use the Mittag-Leffler condition, which is
obviously satisfied when all transition maps are onto. Then, we need to
know that the canonical map
$K_{j+1}(\prod_{i=1}^n \Gr_{d'_i,r'_i})
\to K_{j+1}(\prod_{i=1}^n \Gr_{d_i,r_i}))$ is onto
whenever $d_i\leq d'_i$ and $r_i\leq r'_i$.

An $S$-scheme $X$ is cellular if there exists a sequence of closed
subschemes $\emptyset=Z_0\subset Z_1\subset \dots \subset Z_k=X$ of $S$
such that $Z_i-Z_{i-1}$ is isomorphic to an affine space $\AA^d$ over $S$
for $1\leq i\leq k$.

It is well known that Grassmann varieties are cellular (see \cite{ehresmann})
and it is easy to prove the following formulas:
\begin{itemize}
\item if $X$ is a smooth cellular $S$-scheme, then for any $j\in\NN$,
$K_j(S)\otimes_{K_0(S)} K_0(X)\isomto K_j(X)$;
\item if $X$ is a smooth cellular $S$-scheme, $T$ a regular scheme
and $T\to S$ a morphism, then $K_0(T)\otimes_{K_0(S)}
K_0(X)\isomto K_0(T\times_S X)$;
\item if $X$ and $Y$ are smooth cellular $S$-schemes, then
$K_0(X)\otimes_{K_0(S)} K_0(Y)\isomto K_0(X\times_S Y)$.
\end{itemize}

We see that we only have to prove that $K_0(\Gr_{d',r'})\to
K_0(\Gr_{d,r})$ is onto whenever $d\leq d'$ and $r\leq r'$. One may also
assume that $S=\SpecZ$. Then, for any tuple $(d,r)\in\NN^2$,
$K_0(\Gr_{d,r})$ is generated as a $\lambda$-ring by the class
$u_{d,r}=[\mathcal M'_{d,r}]-d$ (see \cite[VI~4.6]{SGA6}).
With the notations above, the lemma follows from the obvious fact that 
the inverse image of $u_{d',r'}$ by the inclusion
$\Gr_{d,r}\to\Gr_{d',r'}$ is $u_{d,r}$.

\begin{remark}\label{remark-morphisms-pic-to-k-0}
The particular case $d=d'=1$ in the proof shows that the direct
system $(\mathbf{P}^n)_{n\in\mathbf{N}}$ does not unveil phantoms
in the objects $\R\Omega^j(\mathbf{Z}\times\Gr)$. This gives an
interpretation of morphisms $\mathbf{P}^\infty\to \mathbf{Z}\times\Gr$ in
$\Ho$ as natural transformation $\Pic(-)\to K_0(-)$ in $\Sm^\opp\Sets$.
\end{remark}

To finish the proof of theorem~\ref{theorem-computation-k-of-z-gr}, we have
to determine the structure of the ring
$R=\underset{(d,r)\in\NN^2}{\lim}\;K_0(\Gr_{d,r})$. If we fix $d$, we know from
\cite[VI~4.10]{SGA6} that $\underset{r\in\NN}{\lim}\;K_0(\Gr_{d,r})\simeq
K_0(S)[[\tilde{\gamma}^1,\dots,\tilde{\gamma}^d]]$ where
$\tilde{\gamma}^i$ is given by the compatible family $\gamma^i(u_{d,r})$.
Then, $R$ identifies to
$\underset{d\in\mathbf{N}}{\lim}\;K_0(S)[[\tilde{\gamma}^1,\dots,
\tilde{\gamma}^d]]$. One can easily see that the induced transition maps
\[K_0(S)[[\tilde{\gamma}^1,\dots, \tilde{\gamma}^d,\tilde{\gamma}^{d+1}]]
\to K_0(S)[[\tilde{\gamma}^1,\dots, \tilde{\gamma}^d]]\] are obtained by
making $\tilde{\gamma}^{d+1}$ vanish. It proves that $R$ identifies to the
ring of formal power series with an infinite number of variables
$\tilde{\gamma}^1,\tilde{\gamma}^2,\dots$ and coefficient ring $K_0(S)$.

\section{Algebraic structures}
\label{section-algebraic-structures}

We shall see that the previous results show that the algebraic structures
on the sets $K_0(X)$, $X\in\Sm$ uniquely refine to structures of the same
type on $\ZZ\times\Gr$ in the category $\Ho$. Thus, $\ZZ\times\Gr$ shall be
endowed with the structure of a special $\lambda$-ring with duality in
$\Ho$. In this section, we shall use similar notions to those appearing in
\cite{Cohn}.

\subsection{Abstract operators, formulas, algebraic structures}

\begin{definition}
We define a language $\mathcal L$ as the datum of a family of elements
$(l_i)_{i\in I}$ called abstract operators, where each of these operators
 is equipped
with its arity $n_i\in\NN$.
\end{definition}

\begin{definition}
A formula of the language $\mathcal L=(l_i,n_i)_{i\in I}$
involving variables $(x_v)_{v\in V}$
($V$ is assumed to be finite) is the set of expressions inductively built
from the following rules:
\begin{itemize}
\item for any $v\in V$, $x_v$ is a formula;
\item for any $i\in I$, if $F_1,\dots,f_{n_i}$ are formulas, then
$l_i(F_1,\dots,F_{n_i})$ is a formula.
\end{itemize}
\end{definition}

\begin{definition}
An abstract algebraic structure is the datum of a language $\mathcal L$ and
of a family of pairs $(A_r,B_r)_{r\in R}$ of formulas of $\mathcal L$
involving variables in some finite set $V_r$. These pairs are called
``relations'' and shall be denoted $A_r=B_r$.
\end{definition}

\begin{example}\label{example-group-structure}
The abstract algebraic structure of group is defined as follows. The
language $\mathcal L$ is made of a $0$-ary operator $e$ (we may say that $e$ is
a constant), a binary operator $\mu$ and an unary operator $i$. The relations are:
\begin{itemize}
\item $\mu(x,\mu(y,z))=\mu(\mu(x,y),z))$ ;
\item $\mu(e,x)=x$ ;
\item $\mu(x,e)=x$ ;
\item $\mu(x,i(x))=e$ ;
\item $\mu(i(x),x))=e$.
\end{itemize}
Each of these relations involves
a subset of $\{x,y,z\}$ as set of variables.
\end{example}

\subsection{Algebraic structures on objects}

\begin{definition}
Let $\mathcal L=(l_i,n_i)_{i\in I}$ be a language.
An $\mathcal L$-object consists of an object $X$
of a category $\mathcal C$ such that all finite products
$X^n$ exist and of a family of morphisms $X^{n_i}\to X$ denoted $l_i$, for all $i\in I$.
\end{definition}

A morphism of $\mathcal L$-objects $X\to Y$ in a category $\mathcal C$ is a
morphism $F\colon X\to Y$ in $\mathcal C$ such that for any $i\in I$, the
obvious diagram commutes:
\[
\xymatrix{
\ar[d]_{(F,\dots,F)}X^{n_i}\ar[r]^{l_i} & X\ar[d]^F \\
Y^{n_i}\ar[r]^{l_i} & Y}
\]

If $X$ is an $\mathcal L$-object, then one can inductively define a morphism
$F\colon X^V\to X$ for any formula $F$ of $\mathcal L$ involving a finite set
of variables $V$.

\begin{definition}
Let $\mathcal S=(\mathcal L,(A_r=B_r)_{r\in R})$ be an abstract algebraic
structure. An object equipped with an $\mathcal S$-structure is an
$\mathcal L$-object $X$ in some category $\mathcal C$
such that for any $r\in R$, the morphisms
$X^{V_r}\to X$ defined by $A_r$ and $B_r$ are equal. We may also
say that $X$ is an $\mathcal S$-object or that $X$ is a model of $\mathcal
S$ in the category $\mathcal C$.
\end{definition}

We may define the category of $\mathcal S$-objects as a full subcategory
of the category of $\mathcal L$-objects.

\begin{proposition}\label{proposition-lifting-algebraic-structures}
Let $\mathcal S$ be an abstract algebraic structure. Let $F\colon \mathcal
C\to \mathcal D$ be a functor. We assume that finite products exist in
$\mathcal C$ and that $F$ commutes with these products. If $X$ is an
$\mathcal S$-object in $\mathcal C$, then $FX$ has a natural structure of
an $\mathcal S$-object in $\mathcal D$.

Conversely, if the canonical map $\Hom_{\mathcal C}(X^n,X)\to
\Hom_{\mathcal D}(F(X^n),FX)$ is a bijection for any $n\in\NN$ and some
object $X$ of $\mathcal C$, then an $\mathcal S$-structure on $FX$ uniquely
arises from an $\mathcal S$-structure on $X$.

Furthermore, let $X$ and $Y$ be two $\mathcal S$-objects. We assume that
for any $n\in\NN$, the map $\Hom_{\mathcal C}(X^n,Y)\to \Hom_{\mathcal
D}(F(X^n),FY)$ is a bijection. Let $f\colon X\to Y$ be a morphism in
$\mathcal C$. Then, $f$ is a morphism of $\mathcal S$-objects in $\mathcal
C$ if and only if $Ff\colon FX\to FY$ is a morphism of $\mathcal S$-objects
in $\mathcal D$.
\end{proposition}

This is a pliantly true.

\subsection{Structures on $\ZZ\times\Gr$}

The example~\ref{example-group-structure} shows that there is an obvious
abstract algebraic structure whose models in the category of sets are
groups. The
same applies to commutative rings (with unit): the underlying language of
the corresponding abstract algebraic structure involves the $0$-ary
operators $0$ and $1$, the unary operator $-$ and the binary operators $+$
and $\times$. Following \cite[RRR~I~1]{SGA6}, if we add a family of unary
operators $(\lambda^n)_{n\in\NN}$, we can define the abstract algebraic
structures of $\lambda$-rings and of special $\lambda$-rings. One may also
introduce the abstract algebraic structure of special $\lambda$-rings with
duality: we add an unary duality operator that should be an involution
commuting with the other operators.

\begin{theorem}
\label{theorem-z-gr-special-lambda-ring}
Let $S$ be a regular scheme. In the category $\Ho$, there exists a unique
structure of a special $\lambda$-ring with duality on the object
$\ZZ\times\Gr$ such that the corresponding induced structures of
$\lambda$-rings with duality on $K_0(X)$ for
all $X\in\Sm$ are the usual ones.
\end{theorem}

For any $X\in\Sm$, the set $K_0(X)$ is endowed with the structure of
a special $\lambda$-ring with duality \cite[VI~3.2]{SGA6}. All these
structures are compatible with inverse image maps $f^\star\colon K_0(X)\to
K_0(Y)$ for morphisms $f\colon Y\to X$. This shows that, as a presheaf of
sets on $\Sm$, $K_0(-)=\pi_0(\mathbf{Z}\times \Gr)$ is endowed with the
structure of a special $\lambda$-ring with duality.
Proposition~\ref{proposition-lifting-algebraic-structures} and
theorem~\ref{theorem-operations-k-theory-several-operands} shows that it
lifts to a unique structure of a special $\lambda$-ring with duality on
$\ZZ\times \Gr$ in $\Ho$.

\bigskip

\begin{proposition}
Let $f\colon Y\to X$ be a morphism of regular schemes. Let
$\ZZ\times\Gr_X\in \Ho[X]$
(resp. $\ZZ\times\Gr_Y\in \Ho[Y]$) be the special $\lambda$-rings with duality
defined in theorem~\ref{theorem-z-gr-special-lambda-ring}. The
structures on $\ZZ\times\Gr_X$ induce a structure of special
$\lambda$-rings with duality on $\L f^\star (\ZZ\times\Gr_X)$. Then, the
obvious isomorphism $\L f^\star (\ZZ\times\Gr_X)\simeq \ZZ\times\Gr_Y$ in
$\Ho[Y]$ is an isomorphism of special $\lambda$-rings with duality.
\end{proposition}

We can use the construction of
proposition~\ref{proposition-lifting-algebraic-structures} because the
functor $\L f^\star\colon \Ho[X]\to\Ho[Y]$ (see \cite[page~108]{MV})
commutes with finite products. Using
theorem~\ref{theorem-operations-k-theory-several-operands}, it suffices to
compare the two induced special $\lambda$-rings with duality structures on
the presheaf $K_0(-)$ on $\Sm[Y]$. If $f$ is smooth, one may argue by
saying that the structures on $\pi_0\L f^\star (\ZZ\times\Gr_X)$ are
obtained from those on $\pi_0 (\ZZ\times \Gr_X)$ by applying the
``restriction'' functor $\Presh[X]\to\Presh[Y]$ obtained by composition
with the ``forgetful'' functor $\Sm[Y]\to\Sm[X]$. In the general case, we
may observe that it suffices to check that the two special $\lambda$-rings
structures considered on $K_0(-)$ in $\Presh[Y]$ agree on the ``universal''
elements $u_{d,r}+n\in K_0(\Gr_{d,r,Y})$ (see
remark~\ref{remark-explanation-generation-up-to-t}) and this follows from
the fact that the presheaves $K_0(-)$ on $\Sm[X]$ or $\Sm[Y]$ come from a
presheaf of special $\lambda$-rings with duality on the category of all
regular schemes.

\begin{remark}
Similar arguments can be used to prove that, through the interpretation of
operations as formal power series 
(see theorem~\ref{theorem-computation-k-of-z-gr}), the map $\L
f^\star\colon
\End_{\Ho[X]}(\ZZ\times\Gr_X)\to\End_{\Ho[Y]}(\ZZ\times\Gr_Y)$ corresponds
to the extension of scalars of formal power series along the morphism
$f^\star\colon K_0(X)\to K_0(Y)$.
\end{remark}

\subsection{Structures on higher $K$-groups}
\label{subsection-structures-higher-k-groups}

Let $S$ be a regular scheme. We have constructed structures on
$\ZZ\times\Gr$ in $\Ho$.
For any $\mathcal X\in \Ho$, they induce
structures on the set $\Hom_{\Ho}(\mathcal X,\ZZ\times\Gr)$, which we
denote $K_0(\mathcal X)$. As a result, these sets $K_0(\mathcal X)$ are
special $\lambda$-rings with duality. To extend some structures to the
higher $K$-groups $K_n(\mathcal X)=\Hom_{\Hopt}(S^n\wedge \mathcal
X_+,\ZZ\times\Gr)$, one has to refine some morphisms in $\Ho$ to morphisms
in $\Hopt$.

Theorem~\ref{theorem-endomorphisms-z-gr} and the subsequent comments shows
that the families of operations $(\Psi^k)_{k\in\mathbf{Z}}$,
$(\lambda^n)_{n\in\NN}$ and $(\gamma^n)_{n\in\NNstar}$ and more generally
all operations $\tau\colon K_0(-)\to K_0(-)$ such that $\tau(0)=0$ naturally
act on these sets $K_n(\mathcal X)$. Moreover, relations known at the
level of $K_0$ implies similar relations on all the $K$-groups: for
instance, the formula $\Psi^k\circ \Psi^{k'}=\Psi^{kk'}$ is satisfied by
the corresponding operations on $K_\star(\mathcal X)$.

This also applies to operations involving several operands like $+$ and
$\times$. The commutative group structure on $\ZZ\times\Gr$ in $\Ho$ comes
from a commutative group structure on $\ZZ\times\Gr$ in $\Hopt$. Using this
$H$-group structure, we obtain abelian group structures on the sets
$K_n(\mathcal X)$ for all $n\in\NN$. Using the argument of
\cite[page~74]{morel-asterisque},
the product law $\times\colon (\ZZ\times\Gr)^2\to
\ZZ\times\Gr$ in $\Ho$ can easily be refined to a pairing $\mu\colon
(\ZZ\times\Gr)\wedge(\ZZ\times\Gr)\to\ZZ\times\Gr$, which induce pairings
$K_i(\mathcal X)\times K_j(\mathcal Y)\to K_{i+j}(\mathcal X\times\mathcal
Y)$ for $\mathcal X$ and $\mathcal Y$ in $\Ho$. Using this construction in
the case $\mathcal Y=\mathcal X$ and the diagonal morphism $\mathcal
X\to\mathcal X\times\mathcal X$, we get a product law on the graded abelian
group $K_\star(\mathcal X)$. It formally follows from the commutative ring
structure on $\ZZ\times\Gr$ in $\Ho$ that with these definitions,
$K_\star(\mathcal X)$ is a graded commutative ring. One can easily check
compatibilities between the $\lambda$-operations and the product. For
instance, if $k\in\mathbf{Z}$, the fact that $\Psi^k$ is an endomorphism of
the ring $\ZZ\times\Gr$ in $\Ho$ shows that the operation $\Psi^k\colon
K_\star(\mathcal X)\to K_\star(\mathcal X)$ is an endomorphism of graded
rings.

\begin{example}
The following confusing example should warn the reader against
misinterpretations of the previous results. Let $\tau\colon K_0(-)\to
K_0(-)$ be the operation defined by $\tau(x)=x^2$ for any $x\in K_0(X)$ and
$X\in \Sm$. This operation satisfies $\tau(0)=0$; then it induces
maps $\tau\colon K_n(X)\to K_n(X)$ for all $n\in\NN$ and $X\in\Sm$.
However, this operation on higher $K$-groups is unrelated to the squaring
map $K_n(X)\to K_{2n}(X)$ unless $n=0$. Indeed, a simple computation using
the splitting principle shows that $\tau=\Psi^2+2\lambda^2$. To the latter,
we associated maps $K_n(X)\to K_n(X)$ rather than maps $K_n(X)\to K_{2n}(X)$.
\end{example}

\section{Comparison with previous constructions}
\label{section-comparison}
\subsection{Models of algebraic $K$-theory}
\label{subsection-models}

\begin{definition}
Let $S$ be a regular scheme. A candidate model of algebraic $K$-theory
(over $S$) is
an object $\mathcal K\in\Hopt$ equipped \
with a morphism $\alpha_{\mathcal K}\colon K_0(-)\to
\pi_0 \mathcal K$ of presheaves of pointed sets on $\Sm$.
We say that $(\mathcal K,\alpha_{\mathcal K})$
is strict if $\alpha_{\mathcal K}$ is an isomorphism.
\end{definition}

For such an object $\mathcal K$, $\mathcal X\in\Ho$ and $n\in\mathbf{N}$,
we define $K_n^{\mathcal K}(\mathcal X)$ to be the set of morphisms
$\Hom_{\Hopt}(S^n\wedge \mathcal X_+,{\mathcal K})$.

A morphism of candidate models $(\mathcal K,\alpha_{\mathcal K})\to
(\mathcal K',\alpha_{\mathcal K'})$ is the datum of a morphism $f\colon
\mathcal K\to \mathcal K'$ in $\Hopt$ such that $\alpha_{\mathcal
K'}=\pi_0(f)\circ \alpha_{\mathcal K}$.

\begin{proposition}\label{proposition-candidate-models}
Candidate models of algebraic $K$-theory can be associated to the following
definitions of algebraic $K$-theory :
\begin{itemize}
\item Quillen's $Q$-construction \cite[7.1]{quillen-lnm-341};
\item Waldhausen's \cite[\S1.9]{waldhausen};
\item Thomason-Trobaugh's \cite[3.5.3]{thomason-trobaugh}.
\end{itemize}
\end{proposition}

For each of these constructions, there is a well-defined presheaf $\mathcal K$
of pointed simplicial sets of $\Sm$ such that the corresponding $K$-groups are
the homotopy groups of the spaces $\mathcal K(X)$ for all $X\in\Sm$. This
presheaf $\mathcal K$ defines an object in $\Hopt$ and there are canonical maps
for all $X\in\Sm$ (see definition~\ref{definition-pi-zero}) : \[\pi_0 (\mathcal
K(X))\to \Hom_{\Hopt}(X,K)=(\pi_0\mathcal K)(X)\;\text{.}\] For any of these
definitions of algebraic $K$-theory, in degree zero, $\pi_0(\mathcal K(X))$ is
identified to the Grothendieck group $K_0(X)$ of the exact category of vector
bundles on $X$. Then, we get the expected map $\alpha_{\mathcal K}\colon
K_0(-)\to \pi_0\mathcal K$ in $\Sm^\opp\Sets_\bullet$.

\bigskip

Thanks to theorem~\ref{theorem-morel-voevodsky}, the object
$\mathbf{Z}\times\Gr$ is endowed with a structure of a (strict) candidate
model of algebraic $K$-theory. The map $\alpha_{\mathbf{Z}\times\Gr}\colon
K_0(-)\to \pi_0(\mathbf{Z}\times\Gr)$ has the (characteristic) property
that the class $u_{d,r}+n\in K_0(\Gr_{d,r})$ is mapped to the homotopy
class of the inclusion $\Gr_{d,r}\subset \{n\}\times
\Gr\subset\mathbf{Z}\times\Gr$.

The following proposition shows that this model
$(\mathbf{Z}\times\Gr,\alpha_{\mathbf{Z}\times\Gr})$ plays an almost
universal role:

\begin{proposition}
Let $S$ be a regular scheme. Let $(\mathcal K,\alpha_{\mathcal K})$ \
be a candidate
model of algebraic $K$-theory over $S$. Then, there exists a morphism
$(\mathbf{Z}\times\Gr,\alpha_{\mathbf{Z}\times\Gr})\to(\mathcal
K,\alpha_{\mathcal K})$ of candidate models of algebraic $K$-theory.
If this morphism is an isomorphism, then it is unique and we shall say that
$(\mathcal K,\alpha_{\mathcal K})$ is a genuine model of algebraic
$K$-theory.
\end{proposition}

May $\mathcal K$ not be an $H$-group, the surjectivity part of the Milnor
exact sequence stated in lemma~\ref{lemma-milnor-exact-sequence} is still true.
Then, there exists a morphism $f\colon \mathbf{Z}\times\Gr\to \mathcal H$
in $\Hopt$ such that the morphism of presheaves $\alpha_K$ and
$\pi_0(f)\circ \alpha_{\mathbf{Z}\times\Gr}$ in
$\Hom_{\Sm^\opp\Sets}(K_0(-),\pi_0\mathcal K)$ coincide on the universal
classes $u_{d,r}+n\in K_0(\Gr_{d,r})$. Then
remark~\ref{remark-explanation-generation-up-to-t} implies that that they
are equal which proves that $f$ is a morphism of candidate models of
algebraic $K$-theory.

If $f$ is an isomorphism, then we may replace $(\mathcal K,\alpha_{\mathcal
K})$ by $(\mathbf{Z}\times\Gr,\alpha_{\mathbf{Z}\times\Gr})$ and the
uniqueness of $f$ means that there exists a unique endomorphism of
$\mathbf{Z}\times\Gr$ which induces the identity on
$\pi_0(\mathbf{Z}\times\Gr)=K_0(-)$, which is known thanks to
theorem~\ref{theorem-endomorphisms-z-gr}.

\begin{corollary}
Let $S$ be a regular scheme.
If $(\mathcal K,\alpha_{\mathcal K})$ and
$(\mathcal K',\alpha_{\mathcal K'})$ are two genuine models of algebraic
$K$-theory, they are canonically isomorphic and the associated
$K$-groups are also canonically isomorphic for all $\mathcal X\in \Ho$ and
$n\in\mathbf{N}$ :
\[K_n^{\mathcal K}(\mathcal X)\simeq K_n^{\mathcal K'}(\mathcal X)\;\text{.}\]
\end{corollary}

It follows from the fact that both genuine models are canonically
isomorphic to $(\mathbf{Z}\times\Gr,\alpha_{\mathbf{Z}\times\Gr})$.

\bigskip

\begin{proposition}\label{proposition-genuine-models}
Let $S$ be a regular scheme.
The candidate models defined in
proposition~\ref{proposition-candidate-models} are genuine models of
algebraic $K$-theory.
\end{proposition}

The proofs of the comparison theorems between Quillen's, Waldhausen's and
Thomason-Trobaugh's constructions (\cite[1.9]{waldhausen} and
\cite[proposition~~3.10]{thomason-trobaugh}) are functorial enough to imply
that the three corresponding presheaves of pointed simplicial sets induce
isomorphic objects in the pointed homotopy category of the site $\Sm_\Nis$.
Moreover, these objects satisfy the Nisnevich descent property
\cite[theorem~~10.8]{thomason-trobaugh} and the homotopy invariance of
algebraic $K$-theory for regular schemes \cite[\S6]{quillen-lnm-341}
shows that they are $\mathbf{A}^1$-local. As a result, if $\mathcal K$ is
one of these presheaves of pointed simplicial sets, the obvious maps $\pi_n
(\mathcal K(X))\to \Hom_{\Hopt}(S^n\wedge X_+,\mathcal K)$ are bijections
for all $X\in\Sm$.
In particular, the map $\alpha_{\mathcal K}\colon K_0(-)\to \pi_0\mathcal
K$, which is part of the datum of a candidate model, is an isomorphism.
These candidate models are strict ones. Then, the proposition follows from
the fact that the object $\mathcal K$ associated to Quillen's
$Q$-construction is isomorphic to $\mathbf{Z}\times\Gr$, which is implicit
in the proof of theorem~\ref{theorem-morel-voevodsky}.

\subsection{Products}

\begin{proposition}\label{proposition-comparison-waldhausen}
Let $S$ be a regular scheme. For all $X\in\Sm$, $(i,j)\in\mathbf{N}^2$, the
pairing $K_i(X)\times K_j(X)\to K_{i+j}(X)$ defined in
subsection~\ref{subsection-structures-higher-k-groups} is the same as the
one defined by Waldhausen \cite{waldhausen}.
\end{proposition}

First, thanks to the results of subsection~\ref{subsection-models}, it
truly makes sense to say that these pairings coincide as the different
flavours of models of algebraic $K$-theory give canonically isomorphic
groups. Then, as Waldhausen's product on $K_\star(X)$ obviously extends the
standard one on $K_0(X)$,
theorem~\ref{theorem-operations-k-theory-several-operands} shows that we
only need to observe that Waldhausen's pairing is functorial enough to be
defined at the level of presheaves of pointed simplicial sets on $\Sm$ and
thus induces a morphism $\mathcal K\times\mathcal K\to \mathcal K$ in
$\Hopt$ where $\mathcal K$ is the model of algebraic $K$-theory associated
to Waldhausen's definition.

\begin{remark}
Using similar arguments, one may prove that the pairing $K_i(X)\times
K_j(X)\to K_{i+j}(X)$ coincides with the one defined by Quillen (only for
$i=0$ or $j=0$). For $i\neq 0$, $j\neq 0$ and $X$ affine, one may also compare
them with the product defined by Loday using the $+$-construction \cite{loday};
the arguments would be similar to the arguments in
subsection~\ref{subsection-one-operand-operations} below.
In particular, Waldhausen's pairing coincide with
those defined by Quillen and Loday. This comparison was already known
(see \cite{weibel-products}).
\end{remark}

\subsection{Operations involving one operand}
\label{subsection-one-operand-operations}

In his article \cite{soule}, Soul\'e defined an action of
$\Rep_{\mathbf{Z}}\GL=\underset{d\in\mathbf{N}}{\lim}\;\Rep_{\mathbf{Z}}\GL_d$
on the higher algebraic $K$-theory of schemes, where
$\Rep_{\mathbf{Z}}\GL_d$ is the Grothendieck group defined by Serre
\cite{serre-groupes-de-grothendieck}. If we fix a regular base scheme $S$,
theorem~\ref{theorem-main} introduces such an action on $K$-theory of
smooth $S$-schemes for elements $\tau\in\End_{\Presh_\bullet}(K_0(-))$. As
we would like to state a compatibility between these two constructions, we
shall introduce a common input for both of them.

\begin{definition}
Let $d\in\mathbf{N}$. We let $\Univ d$ be the universal special
$\lambda$-ring equipped with an element $\ID d$ satisfying the following
conditions:
\begin{itemize}
\item[(i)] $\lambda^d(\ID d)$ is invertible;
\item[(ii)] $\lambda^k(\ID d)$ vanishes for $k\geq d+1$.
\end{itemize}
The special $\lambda$-ring $\Univ \infty$ is the projective limit of the
system $(\Univ d)_{d\in\mathbf{N}}$ where the transition map $\Univ
{d+1}\to \Univ d$ maps $\ID {d+1}$ to $\ID d+1$.
\end{definition}

Obviously, for any $d$, there is a canonical morphism of special
$\lambda$-rings $\Univ d\to \Rep_{\mathbf{Z}}\GL_d$ that maps $\ID d$ to the
class of the tautological representation $\id\colon \GL_d\to\GL_d$ of rank
$d$ of the group scheme $\GL_d$. Serre's computation
\cite[\S3.8]{serre-groupes-de-grothendieck} shows that this sequence of 
morphisms consists of isomorphisms.
Then, the canonical morphism $\Univ \infty\to
\Rep_{\mathbf{Z}}\GL$ is an isomorphism.

We may also use the universal properties of the special $\lambda$-rings
$\Univ d$ to define a morphism of special $\lambda$-rings $\Univ \infty\to
K_0(\Gr)=\Hom_{\Ho}(\Gr,\mathbf{Z}\times\Gr)$. It is induced by the
morphisms
\[\Univ d \to \Hom_{\Ho}(\Gr_{d,\infty},\mathbf{Z}\times\Gr)\simeq
\underset{r\in\mathbf{N}}{\lim}\;K_0(\Gr_{d,r})\]
sending $\ID d$ to the compatible family of classes $([\mathcal
M'_{d,r}])_{r\in\mathbf{N}}$ (see the proof of
lemma~\ref{lemma-generation-up-to-t-of-z-gr} for this notation).

We let $(\Univ \infty)_0$ and $(\Rep_{\mathbf{Z}}\GL)_0$ 
be the kernel of
the rank morphism from these groups to $\mathbf{Z}$. Similarly, we denote
$\tilde{K}_0(\Gr)$ the kernel of the restriction to the base-point
$K_0(\Gr)\to K_0(S)$. The comparison
theorem announced above is the following:

\begin{theorem}
\label{theorem-compatibility-soule}
Let $S$ be a regular scheme. For any $n\geq 1$, the following diagram
commutes:
\[
\xymatrix{& (\Univ \infty)_0\ar[dl]_{\sim}\ar[dr] &  \\
(\Rep_{\mathbf{Z}}\GL)_0\ar[rd] & & 
\ar[ld]\tilde{K}_0(\Gr)\\
&\End_{\Sm^\opp\Ab}(K_n(-))\;\text{,}&}\]
where the two upper maps are the ones mentioned above, the lower-left
one is the one defined by Soul\'e and the lower-right one arises from
theorem~\ref{theorem-morel-voevodsky}.
(Thanks to previous results, this
lower-right map can be interpreted as the canonical map
$\End_{\Presh_\bullet}(\tilde{K}_0(-))\to
\End_{\Sm^\opp\Ab}(K_n(-))$.)
\end{theorem}

The strategy of the proof consists in the construction of an horizontal map
$(\Rep_{\mathbf{Z}}\GL)_0\to \tilde{K}_0(\Gr)$ which makes both upper and
lower triangles commute. This map is induced by a morphism
$\Rep_{\mathbf{Z}}\GL\to K_0(\Gr)$ and is a particular case of a more
general construction:

\begin{proposition}
\label{proposition-representation-ring}
Let $G$ be a smooth group scheme over $\SpecZ$. We let $\Rep_{\mathbf{Z}}G$
be the Grothendieck group of finitely generated free $\mathbf{Z}$-modules
endowed with a linear action of $G$ (see
\cite[\S2.3]{serre-groupes-de-grothendieck}). We let
$\mathbf{B}G\in\Hopt$ be the classifying space of $G$ (where $G$ is
considered as a sheaf of groups on $\Sm_\Nis$). Let $\rho\colon G\to
\GL(M)$ be a free finitely generated $\mathbf{Z}$-module endowed with a
linear action of $G$ (we shall say that $M$ is a representation of $G$). The
choice of a $\mathbf{Z}$-basis of $M$ identifies $\rho$ with a morphism of
group schemes $G\to \GL_d$ over $\SpecZ$ where $d=\rk M$. We let
$[\rho]\in \Hom_{\Ho}(\mathbf{B}G,\mathbf{Z}\times\Gr)$ be the morphism
obtained from $\mathbf{B}\rho\colon \mathbf{B}G\to \mathbf{B}\GL_d$ by
composing with the canonical morphism $\mathbf{B}\GL_d\simeq
\Gr_{d,\infty}\simeq \left\{d\right\}\times\Gr_{d,\infty}\subset
\mathbf{Z}\times\Gr$.
Then, this assignment $\rho\longmapsto [\rho]$ does not depend on the
the choice of $\mathbf{Z}$-bases and induces a morphism of special
$\lambda$-rings with duality
$\Rep_{\mathbf{Z}}G\to\Hom_{\Ho}(\mathbf{B}G,\mathbf{Z}\times\Gr)=
K_0(\mathbf{B}G)$.
\end{proposition}

The choice of two different $\mathbf{Z}$-bases of a representation $M$ of
$G$ would lead to morphisms $G\to \GL_d$ which would differ by an inner
automorphism of $\GL_d$ (induced by an element of $\GL_d(\mathbf{Z})$): the
associated morphisms $\mathbf{B}G\to\mathbf{B}\GL_d$ are equal in $\Ho$
(and also in $\Hopt$ after composition with
$\mathbf{B}\GL_d\to\mathbf{B}\GL_\infty$ because $\mathbf{B}\GL_\infty$ is
an $H$-group).

To prove that $\rho\longmapsto [\rho]$ induces a morphism at the level of
the Grothendieck group of representations of $G$, we use the following two
lemmas:

\begin{lemma}
\label{lemma-sum-on-bgl}
Let $+\colon
\mathbf{B}\GL_\infty\times\mathbf{B}\GL_\infty\to\mathbf{B}\GL_\infty$ be
the $H$-group structure coming from the usual group structure on
$\tilde{K}_0(-)$ (see remark~\ref{remark-variant-gr}).
For any $d,d'\geq 0$, the following diagram commutes in $\Hopt$:
\[
\xymatrix{
\ar[d]
\mathbf{B}\GL_d\times\mathbf{B}\GL_{d'}\ar[r]^-{\mathbf{B}\oplus} &
\mathbf{B}\GL_{d+d'}\ar[d]\\
\mathbf{B}\GL_\infty\times\mathbf{B}\GL_\infty\ar[r]^-{+} &
\mathbf{B}\GL_\infty\;\text{,}}
\]
where the vertical morphisms are the obvious ones and the upper one is the
morphism $\mathbf{B}\oplus$ deduced from the ``direct sum'' morphism
$\oplus\colon \GL_d\times\GL_{d'}\to\GL_{d+d'}$.
\end{lemma}

The correspondence between $\GL_d$-torsors on schemes and rank-$d$ vector
bundles provides a functorial map $H^1(X,\GL_d)\to \tilde{K}_0(X)$ (we
substract the rank in $K_0(X)$ so as to get elements in $\tilde{K}_0(X)$).
An obvious verification leads to the following commutative square which
states a compatibility between this correspondence, the sum in
$\tilde{K}_0(X)$ and the map induced on cohomology by the morphism
$\oplus\colon \GL_d\times\GL_{d'}\to\GL_{d+d'}$:

\[
\xymatrix{H^1(X,\GL_d)\times H^1(X,\GL_{d'})\ar[d] \ar[r]^-{\oplus_\star} &
H^1(X,\GL_{d+d'}) \ar[d]\\
\tilde{K}_0(X)\times \tilde{K}_0(X)\ar[r]^-{+} & \tilde{K}_0(X)}
\]

The two morphisms we want to compare are in
$\Hom_{\Ho}(\mathbf{B}\GL_d\times\mathbf{B}\GL_{d'},\mathbf{B}\GL_\infty)\simeq
\lim_{(r,r')} \tilde{K}_0(\Gr_{d,r}\times\Gr_{d',r'})$. Then, the lemma
follows from the commutativity mentioned above in the case where $X$ is 
a product of Grassmann varieties and where the torsors corresponds to the
universal vector bundles on these varieties.

\begin{lemma}
Let $0\to\rho'\to \rho\to\rho''\to 0$ be an exact sequence of
representations of $G$. Then, $[\rho]=[\rho'\oplus\rho'']$ in
$\Hom_{\Ho}(\mathbf{B}G,\mathbf{Z}\times\Gr)$.
\end{lemma}

Let $d'=\rk \rho'$, $d=\rk\rho$ and $d''=\rk \rho''$. Using the obvious
functoriality of the constructions with respect to the group $G$, we may
assume that we are in the universal situation where $\rho\colon G \to
\GL_d$ is the inclusion of the subgroup of matrices of the form
$g=\left(\begin{array}{cc} g' & h \\ 0 & g'' \end{array}\right)$ where
$g'\in \GL_{d'}$, $g''\in \GL_{d''}$ and $h$ is an $d'$-by-$d''$ matrix and
where the representations $\rho'$ and $\rho''$ correspond to the obvious
morphisms $G\to \GL_{d'}$ and $G\to \GL_{d''}$.

Let $D=\GL_{d'}\times\GL_{d''}$ be the subgroup of $G$ consisting of
matrices of the previous form such that $h=0$. Obviously, the restriction
of the representations $\rho$ and $\rho'\oplus\rho''$ from $G$ to $D$ are
isomorphic. Then, to finish the proof, it suffices to know that the
restriction map $K_0(\mathbf{B}G)\to K_0(\mathbf{B}D)$ is an injection.
Indeed, this map is a bijection because $D\to G$ is an $\mathbf{A}^1$-weak
equivalence and thus $\mathbf{B}D\to\mathbf{B}G$ is also an
$\mathbf{A}^1$-weak equivalence
(see \cite[Proposition~2.14, page~74]{MV}).

\medskip

We have constructed a morphism of abelian groups $\Rep_{\mathbf{Z}}G\to
K_0(\mathbf{B}G)$. To finish the proof of the proposition, it remains to
show that this is a morphism of special $\lambda$-rings with duality. The
compatibility of the construction with external powers and duality can be
checked in the same way as we did it for direct sums (see
lemma~\ref{lemma-sum-on-bgl}).

\bigskip

To prove theorem~\ref{theorem-compatibility-soule},
we apply proposition~\ref{proposition-representation-ring}
to the cases $G=\GL_d$ for
all $d$. It provides a morphism of special $\lambda$-rings
$\Rep_{\mathbf{Z}}\GL_d\to K_0(\Gr_{d,\infty})$. Taking the projective
limit over all $d$ and considering the rank-$0$ part leads to the expected
morphism $(\Rep_{\mathbf{Z}}\GL)_0\to \tilde{K}_0(\Gr)$. The universal
property of $\Univ d$ and the fact that the morphisms
$\Rep_{\mathbf{Z}}\GL_d\to K_0(\Gr_{d,\infty})$ are morphisms of special
$\lambda$-rings shows that the upper triangle commutes. The fact that the
lower triangle commutes follows easily from the very definition in Soul\'e's
paper \cite{soule}.

\section{Virtual categories}
\label{section-virtual-categories}

Virtual categories were introduced by Deligne in
\cite{deligne-determinant}. They are refinements of $K_0$-groups. More
precisely, if $X\in\Sm$ ($S$ regular), the category $\mathcal V(X)$ is
identified to the fundamental groupoid of $\mathcal K(X)$ where $\mathcal
K$ is some $\mathbf{A}^1$-fibrant genuine model of algebraic $K$-theory.
Any vector bundle $\mathcal E$ on $X$ defines an object $\mathcal E$
of the category $\mathcal V(X)$ whose isomorphism class corresponds to
$[\mathcal E]$ in $K_0(X)$. When we have a short exact sequence $0\to
\mathcal E'\to \mathcal E\to \mathcal E''\to 0$, we not only have an
equality of classes $[\mathcal E]=[\mathcal E'\oplus\mathcal E'']$, which
means that $\mathcal E$ and $\mathcal E'\oplus\mathcal E''$ become
isomorphic in $\mathcal V(X)$ but we have a specific isomorphism $\mathcal
E\simeq \mathcal E'\oplus\mathcal E''$ in this category $\mathcal V(X)$.

\subsection{The Thom spectrum of a virtual bundle}

The construction of this paragraph will be used only in
\S\ref{subsubsection-morphisms-f-star-bgl-to-bgl}. It appears here because
it favours the understanding of virtual categories.

\begin{proposition}
\label{proposition-thom-spectrum}
Let $X$ be a scheme. The construction of the Thom spectrum $\Th_X\mathcal
E$ of a vector bundle $\mathcal E$ on $X$ (see
\cite[Definition~2.16, page~111]{MV}) extends to a functor $\Th_X\colon
\mathcal V(X)\to \SH[X]$.
\end{proposition}

(See also \cite[Th\'eor\`eme 1.5.18]{ayoub-asterisque-i}.)
One may first check that the Thom spectrum of a vector bundle is invertible
for the $\wedge$-product in $\SH[X]$; one is reduced to the case of a trivial
bundle because the invertibility can be checked locally for the Zariski
topology on $X$. Then, using the universal property of $\mathcal V(X)$ as a
Picard category, one has to define an isomorphism $\Th_X\mathcal E'\wedge
\Th_X\mathcal E''\simeq \Th_X\mathcal E$ for any short exact sequence
$0\to\mathcal E'\to\mathcal E\to\mathcal E''\to 0$ of vector bundles. If
the sequence splits, a splitting of it gives such an isomorphism (see
\cite[Proposition~2.17, page 112]{MV}) and 
explicit $\mathbf{A}^1$-homotopies show that it is independant of the
splitting. The general case reduces to this because we can use a torsor
$T\to X$ under a vector bundle such that the inverse image of the sequence
splits over $T$. If Jouanolou's trick is available (see
theorem~\ref{theorem-jouanolou-s-trick}), we may use it; otherwise, as I
learned from Dennis Eriksson, we can always use the scheme which
parametrises the sections of $\mathcal E\to\mathcal E''$: it is a torsor
under the vector bundle $\SheafHom(\mathcal E'',\mathcal E')$. To check the
needed coherence properties, we may split a finite number of short exact
sequences of vector bundles as above; then, it becomes straightforward.

\begin{definition}
\label{definition-thom-spectrum}
Let $f\colon X\to S$ be a smooth morphism between noetherian schemes.
Proposition~\ref{proposition-thom-spectrum} defines a functor $\Th_X\colon
\mathcal V(X)\to \SH[X]$. We also denote $\Th_X\colon \mathcal V(X)\to \SH$
the functor obtained by composition with $\L f_\sharp\colon \SH[X]\to \SH$
(see \cite[Proposition~4.4]{riou-smf}).
\end{definition}

\subsection{Inverting primes on $\mathbf{Z}\times\Gr$}

\begin{definition}
Let $S$ be a regular scheme. Let $a\in\mathbf{N}-\{0\}$. For any
$(d,r)\in\mathbf{N}^2$, we define a morphism $\Gr_{d,r}\to\Gr_{ad,ar}$
which sends an admissible subbundle $\mathcal M'\subset \mathcal O^{n}$
($n=d+r$) of rank $d$ to $\delta_{a,n}(\mathcal M'^{\oplus a})$ where
$\delta_{a,n}\colon (\mathcal O^{n})^{\oplus a}\to \mathcal O^{an}$ is the
isomorphism that sends $(s^1_1,\dots,s^1_n),\dots,(s^a_1,\dots,s^a_n)$ to
$(s^1_1,\dots,s^a_1,\dots,s^1_n,\dots,s^a_n)$. This compatible family of
morphisms induces a morphism $m_a\colon \Gr\to \Gr$ of presheaves of
pointed sets. We also denote
$m_a\colon\mathbf{Z}\times\Gr\to\mathbf{Z}\times\Gr$ the morphism which is
the multiplication by $a$ on $\mathbf{Z}$ and $m_a$ on $\Gr$.
\end{definition}

\begin{lemma}\label{lemma-multiplication-on-z-gr}
Let $a$ and $b$ be two positive natural numbers. Then, the endomorphisms
$m_{ab}$ and $m_a\circ m_b$ of $\mathbf{Z}\times\Gr$ are equal in the
category of presheaves of pointed sets.
\end{lemma}

\begin{definition}
For any $x\in\mathbf{N}-\{0\}$, we set $\frac
1x(\mathbf{Z}\times\Gr)=\mathbf{Z}\times\Gr$. If $y\in\mathbf{N}-\{0\}$
is a multiple of $x$, the endomorphism $m_{y/x}$ of $\mathbf{Z}\times\Gr$
defines a canonical morphism $\frac 1x(\mathbf{Z}\times\Gr)\to\frac
1y(\mathbf{Z}\times\Gr)$.
\end{definition}

Lemma~\ref{lemma-multiplication-on-z-gr} says that this defines a direct
system $(\frac 1x(\mathbf{Z}\times \Gr))_{x\in\mathbf{N}-\{0\}}$ of sheaves
of pointed sets of $\Sm$. It is indexed by $\mathbf{N}-\{0\}$, which is
ordered by divisibility.

\begin{definition}
Let $n$ be a supernatural number (see
\cite[\S{}I.1.3]{cohomologie-galoisienne}). We denote 
$(\mathbf{Z}\times\Gr)[\frac 1 n]$ the colimit of the system
$\frac 1 x (\mathbf{Z}\times\Gr)$ where $x$ varies in the set of positive
natural numbers dividing $n^\infty$.
\end{definition}

\begin{proposition}\label{proposition-endomorphisms-z-gr-one-over-n}
Let $S$ be a regular scheme. Let $i\in\mathbf{N}$. Let $n$ be a
supernatural number. Then, the canonical maps are bijections:
\[\Hom_{\Ho}(\mathbf{Z}\times\Gr,\R\Omega^i(\mathbf{Z}\times\Gr)[{\textstyle
\frac 1 n}])\isomto
\Hom_{\Sm^\opp\Sets}(K_0(-),K_i(-)[{\textstyle\frac 1 n}])\;\text{,}\]
\[\Hom_{\Ho}((\mathbf{Z}\times\Gr)[{\textstyle \frac 1n}],\R\Omega^i(\mathbf{Z}\times\Gr)[{\textstyle
\frac 1 n}])\isomto
\Hom_{\Sm^\opp\Sets}(K_0(-)[{\textstyle \frac 1n}],
K_i(-)[{\textstyle\frac 1 n}])\;\text{.}\]
\end{proposition}

A variant of lemma~\ref{lemma-property-k-for-z-gr} shows that
$\R\Omega^i(\mathbf{Z}\times\Gr)[\frac 1n]$ satisfies property~(K). Hence,
theorem~\ref{theorem-property-k-implies-computation} gives the first
bijection. The second bijection needs additional arguments.
From lemma~\ref{lemma-generation-up-to-t-of-z-gr}, it is easy to show that
$(\mathbf{Z}\times\Gr)[\frac 1n]$ generates
$\pi_0((\mathbf{Z}\times\Gr)[\frac 1n])$ up to $\mathcal T$. The definition
also gives an expression of $(\mathbf{Z}\times\Gr)[\frac 1n]$ as the
colimit of some direct system $(X_i)_{i\in\mathcal I}$ of representable
sheaves, where $\mathcal I$ is an ordered set which has a cofinal sequence.
Then, using proposition~\ref{proposition-key},
we have to show that $(X_i)_{i\in\mathcal I}$ does not unveil
phantoms in $\R\Omega^i(\mathbf{Z}\times\Gr)[\frac 1n]$. Reasoning like in
the proof of lemma~\ref{lemma-property-k-for-z-gr}, it suffices to 
check that for a natural number $a$ dividing a power of $n$, the
morphisms $m_{a,d,r}\colon \Gr_{d,r}\to\Gr_{ad,ar}$ induce surjections
$m_{a,d,r}^\star\colon
K_0(\Gr_{ad,ar})[\frac 1n]\to K_0(\Gr_{d,r})[\frac 1n]$. This is true
because $K_0(\Gr_{d,r})[\frac 1n]$ is generated by $u_{d,r}$ as
a $K_0(S)[\frac 1n]$-$\lambda$-algebra and 
$m_{a,d,r}^\star(u_{ad,ar})=au_{d,r}$.

We leave the variants involving several operands to the reader.

\subsection{Operations on virtual categories}

\begin{definition}
Let $S$ be a regular scheme. Let $n$ be a supernatural number.
We let $\mathcal V(-)[\frac 1n]$ be the
presheaf of groupoids that sends $X\in\Sm$ to the fundamental groupoid
of $\mathcal K[\frac 1n](X)$ where $\mathcal K[\frac 1n]$ is an
$\mathbf{A}^1$-fibrant replacement of $(\mathbf{Z}\times\Gr)[\frac 1n]$.
\end{definition}

\begin{theorem}
\label{theorem-operation-on-virtual-categories}
Let $n$ be an even supernatural number. Let $\tau\colon K_0(-)[\frac 1n]\to
K_0(-)[\frac 1n]$ be a morphism in $\Sm[\SpecZ]^\opp\Sets$. Then, up to a
unique isomorphism, we can define a family of functors
$\tilde{\tau}_X\colon \mathcal V(X)[\frac 1n]\to\mathcal V(X)[\frac
1n]$ for $X\in\Sm[\SpecZ]$ which induces $\tau$ on sets of isomorphisms
classes in $\mathcal V(-)[\frac 1n]$ and such that for any morphism
$f\colon Y\to X$ in $\Sm[\SpecZ]$, it satisfies the equality $f^\star\circ
\tilde{\tau}_X=\tilde{\tau}_Y\circ f^\star$. (Variants involving several
operands are also true.)
\end{theorem}

From proposition~\ref{proposition-endomorphisms-z-gr-one-over-n}, we know
that $\tau$ corresponds to an endomorphism of $\mathcal K[\frac 1n]$ in
$\Ho[\SpecZ]$. As $\mathcal K[\frac 1n]$ is $\mathbf{A}^1$-fibrant, $\tau$
lifts to a morphism $\tilde{\tau}\colon \mathcal K[\frac 1n]\to\mathcal
K[\frac 1n]$. Passing to fundamental groupoids, we get a family of functors
$\tilde{\tau}_X\colon \mathcal V(X)[\frac 1n]\to \mathcal V(X)[\frac 1n]$
for $X\in\Sm[\SpecZ]$.

We let $E=\Sheafhom(\mathcal K[\frac 1n],\mathcal K[\frac 1n])$ be the
simplicial set of endomorphisms of $\mathcal K[\frac 1n]$ (it is given by
the simplicial structure). The morphism $\tilde{\tau}$ corresponds to a
$0$-simplex in $E$. If $\tilde{\tau}'\colon \mathcal K[\frac 1n]\to\mathcal
K[\frac 1n]$ is in the same homotopy class as $\tilde{\tau}$, the choice of
an homotopy (\ie, a path between $\tilde{\tau}$ and $\tilde{\tau}'$ in $E$ 
gives an isomorphism between the associated
families of functors $(\tilde{\tau}_X)$ and
$(\tilde{\tau}'_X)$. The question is whether this isomorphism is uniquely
determined or not. It will be so if there exists a unique homotopy class of
paths $\tilde{\tau}\to\tilde{\tau}'$. As $E$ is an $H$-group, it means that
the connected components of $E$ are simply connected, \ie, $\pi_1E=0$. This
group identifies to $\Hom_{\Ho[\SpecZ]}(\mathcal K[\frac 1n],\Omega\mathcal
K[\frac 1n])$, which identifies to
$\Hom_{\Sm[\SpecZ]^\opp\Sets}(K_0(-)[\frac 1n],K_1(-)[\frac 1n])$.
To prove that this group vanishes, we can use
proposition~\ref{proposition-key} which expresses it as a projective limit
of some groups $K_1(\Gr_{d,r})[\frac 1n]$. The result then follows from the
fact that $K_1(\mathbf{Z})\simeq \mathbf{Z}/2\mathbf{Z}$.

\begin{remark}
In theorem~\ref{theorem-operation-on-virtual-categories}, we may replace
$\Sm[\SpecZ]$ by any small full subcategory $\Reg$ of the category of regular
schemes. Indeed, we may assume that $\SpecZ\in\Reg$ and that for any
$S\in\Reg$, objects in $\Sm$ belong to $\Reg$. Then, we may work in the
$\mathbf{A}^1$-homotopy category $\Ho[\Reg]$ of the site $\Reg_\Nis$
equipped with the interval $\mathbf{A}^1$. Arguments leading to
theorem~\ref{theorem-operation-on-virtual-categories} can be made with the
category $\Ho[\Reg]$ instead of $\Ho[\SpecZ]$. We may also deduce 
results for $\Reg$ from the case of $\Sm[\SpecZ]$ by using the fully
faithful functor $\L p^\star\colon \Ho[\SpecZ]\to\Ho[\Reg]$ associated to
the obvious reasonable continuous map of sites
$p\colon \Reg_\Nis\to \Sm[\SpecZ]_\Nis$.
\end{remark}

\section{Additive and stable results}
\label{section-additive-and-stable-results}
\subsection{The splitting principle}

Now, we shall focus on natural transformations $K_0(-)\to K_0(-)$ which are
compatible with the abelian group structures on $K$-groups, \ie, morphisms
in $\Sm^\opp\Ab$ rather than in $\Presh$. From
theorem~\ref{theorem-operations-k-theory-several-operands} and
proposition~\ref{proposition-lifting-algebraic-structures}, we know that
these additive operations precisely correspond to endomorphisms of
$\mathbf{Z}\times\Gr$ as an $H$-group (\ie, a group object in $\Hopt$).

To compute these additive transformations, we shall use the ``splitting
principle''. We let $\Pic(-)$ be the presheaf of sets on $\Sm$ (for a regular
scheme $S$) that maps $U\in\Sm$ to the Picard group $\Pic(U)$, considered
as a set. We denote $c\colon \Pic(-)\to K_0(-)$ the morphism in $\Presh$
that maps the isomorphism class of a line bundle $\mathcal L$ to the class
$[\mathcal L]$ in the Grothendieck group of vector bundles.

\begin{proposition}
\label{proposition-splitting-principle}
Let $S$ be a regular scheme. For any integer $i$, the map induced by $c$
\[c^\star\colon
\Hom_{\Sm^\opp\Ab}(K_0(-),K_i(-))\to\Hom_{\Sm^\opp\Sets}(\Pic(-),K_i(-))\]
is a bijection. Moreover, the latter group identifies to
\[\Hom_{\Ho}(\mathbf{P}^\infty, \R\Omega^i(\mathbf{Z}\times\Gr))\simeq
\lim_n K_i(\mathbf{P}^n)\simeq K_i(S)[[U]]\;\text{,}\]
where $U=[\mathcal O(1)]-1$ is the obvious compatible family in $\lim_n
K_0(\mathbf{P}^n)$.
\end{proposition}

The injectivity of $c^\star$ follows easily from the ``splitting
principle'': if $\mathcal M$ is a vector bundle of rank $r$ on a scheme
$X\in\Sm$, the complete flag scheme $\mathbf{D}(\mathcal
M)\vers{\pi} X$ is such that $[\pi^\star \mathcal M]$ decomposes in
$K_0(\mathbf{D}(\mathcal M))$ as a sum of the classes of $r$ line bundles
and $\pi^\star\colon K_i(X)\to K_i(\mathbf{D}(\mathcal M))$ is injective.

Proposition~\ref{proposition-key},
lemma~\ref{lemma-generation-up-to-t-of-z-gr},
lemma~\ref{lemma-property-k-for-z-gr} and
remark~\ref{remark-morphisms-pic-to-k-0}
show that we have bijections:
\[ \Hom_{\Sm^\opp\Sets}(\Pic(-),K_i(-))
\simeq \Hom_{\Ho}(\mathbf{P}^\infty, \R\Omega^i(\mathbf{Z}\times\Gr))\simeq
\lim_n K_i(\mathbf{P}^n)\;\text{.}\]
The identification of this group with $K_i(S)[[U]]$ follows from the
computation of the algebraic $K$-theory of projective spaces:
$K_i(\mathbf{P}^n)\simeq K_0(\mathbf{P}^n)\otimes_{K_0(S)}K_i(S)$ and
$K_0(\mathbf{P}^n)\simeq K_0(S)[U]/(U^{n+1})$.

It remains to show that $c^\star$ is surjective. Using the previous
identifications, we rewrite it as a map $c^\star\colon
\Hom_{\Sm^\opp\Ab}(K_0(-),K_i(-))\to K_i(S)[[U]]$. First, we observe that
for any $k\in\mathbf{N}$ and $x\in K_i(S)$, we may denote $x\Psi^k$ the
natural transformation $K_0(-)\to K_i(-)$ that maps $y$ to $x\cdot
\Psi^k(y)$ and see that it satisfies $c^\star(x\Psi^k)=x(1+U)^k$. This
proves that the image of $c^\star$ contains $K_i(S)[U]$. To finish the
proof, we use the following lemma:

\begin{lemma}\label{lemma-convergence-operations}
Let $(\tau_n)_{n\in\mathbf{N}}$ be a sequence of additive natural
transformations $K_0(-)\to K_i(-)$ such that $c^\star(\tau_n)$ converges to
zero in $K_i(S)[[U]]$ for the infinite product topology, where
$K_i(S)$ is endowed with the discrete topology; in other words, we assume
that for each $k\in\mathbf{N}$, the coefficient of $U^k$ in
$c^\star(\tau_n)$ eventually vanishes. Then, for any $X\in\Sm$ and $x\in
K_0(X)$, there exists $N\in\mathbf{N}$ such that for all $n\geq N$,
$\tau_n(x)=0$ and it makes sense to define a natural transformation
$\tau\colon K_0(-)\to K_i(-)$ by the formula $\tau(x)=\sum_{n=0}^\infty
\tau_n(x)$ and we have the equality
$c^\star(\tau)=\sum_{n\in\mathbf{N}} c^\star(\tau_n)$ in $K_i(S)[[U]]$.
\end{lemma}

We have to prove that given $X\in\Sm$ and $x\in K_0(X)$, $\tau_n(x)$
eventually vanishes. The assumption says that it is true for the class
$x=[\mathcal O(1)]$ on $\mathbf{P}^n$ for all $n$. Taking inverse images of
these classes by morphisms $f\colon X\to\mathbf{P}^n$ enables to obtain the
more general case of classes of line bundles generated by their global
sections, \eg, line bundles on affine schemes. Using $\mathcal T$ (see
theorem~\ref{theorem-jouanolou-s-trick}), we get the case of line bundles
on any $X\in\Sm$. Then, the general case follows from the splitting
principle.

\begin{remark}\label{remark-topology-operations}
We may define a topology on $\Hom_{\Sm^\opp\Ab}(K_0(-),K_i(-))$ by
considering the weakest topology for which the evaluation maps at all
elements $x\in K_0(X)$ for all $X\in\Sm$ are continuous, where all groups
$K_i(X)$ are endowed with the discrete topology. The argument of the lemma
shows that the bijection $c^\star\colon
\Hom_{\Sm^\opp\Ab}(K_0(-),K_i(-))\to K_i(S)[[U]]$ is an homeomorphism.
\end{remark}

\begin{remark}\label{remark-law-bigstar}
The composition of endomorphisms endows \[\End_{\Sm^\opp\Ab}(K_0(-))\simeq
K_0(S)[[U]]\] with a structure of a (possibly non-commutative) ring. If this
law on $K_0(S)[[U]]$ is denoted $\bigstar$, one may characterise it by the
fact that it is continuous and that for all $(x,y)\in K_0(S)^2$,
$(k,k')\in\mathbf{N}^2$, $(x(1+U)^k)\bigstar(y(1+U)^{k'})=
(x\Psi^k(y)(1+U)^{kk'}$. More generally, we have a graded ring structure
on $\oplus_{i\in\mathbf{N}}K_i(S)[[U]]$ which comes from the fact that
$K_i(S)[[U]]$ identifies to the group of homomorphisms
$\mathbf{Z}\times\Gr\to \R\Omega^i(\mathbf{Z}\times\Gr)$
of abelian groups inside $\Ho$; the multiplication can be described
similarly as it has been described in degree $0$.
\end{remark}

\begin{remark}\label{remark-splitting-principle}
The surjectivity of the map
\[c^\star\colon\Hom_{\Sm^\opp\Ab}(K_0(-),K_i(-))\to K_i(S)[[U]]\]
may be proved using a different argument. First, we may assume that the
formal power
series $f=\sum_{k\geq 0}a_kU^k$ is such that $a_0=0$, so that we actually
have to construct a natural transformation $\tau\colon \tilde{K}_0(-)\to
K_i(-)$. After the application of Jouanolou's trick and the
splitting principle, an element $x\in \tilde{K}_0(X)$ can be expressed as
$x=\sum_{i=1}^n u_i$ where $u_i=[\mathcal L_i]-1$ for a family of vector
bundles $L_1,\dots,L_n$ for a big enough $n$. Then, $\tau(x)$ should be 
$f(u_1)+\dots+f(u_n)=\sum_{k\geq 1}a_k(u_1^k+\dots+u_n^k)$. Using the
theory of symmetric polynomials and once we have noticed that the elementary
symmetric functions of the $u_1,\dots,u_n$ are the elements
$\gamma^1(x),\dots,\gamma^n(x)$, we get the existence of an element in
$K_i(S)[[\tilde{\gamma}^1,\tilde{\gamma}^2,\dots]]$
whose associated natural transformation
$\tau\colon \tilde{K}_0(-)\to K_i(-)$ (see
theorem~\ref{theorem-computation-k-of-z-gr}) is additive and such that
$\tau(u)=f(u)$ whenever $u=[\mathcal L]-1$
and $\mathcal L$ is a line bundle.
\end{remark}

\begin{remark}
The method of remark~\ref{remark-splitting-principle} can be used for the
study of natural transformation $\tau\colon K_0(-)\to K_0(-)$ which induces
group morphisms $(K_0(X),+)\to (K_0(X)^\times,\cdot)$ for all $X\in\Sm$,
\ie, classes which are multiplicative on short exact sequences (\eg, the
Todd class, which is defined after tensoring with $\mathbf{Q}$). The result
is that for any series $f=\sum_{k\geq 0} a_k U^k$ such that $a_0$ is
invertible in $K_0(S)$, there exists a unique $\tau\colon K_0(-)\to K_0(-)$
as above such that for any line bundle $\mathcal L$, $\tau([\mathcal
L])=f(u)$ where $u=[\mathcal L]-1$. The proof follows the same pattern:
reduce to the case $a_0=1$ and then consider $f(u_1)f(u_2)\dots f(u_n)$
instead of $f(u_1)+f(u_2)+\dots +f(u_n)$.
\end{remark}

\begin{exercise}[Optional]
Assume that $S$ is a regular scheme such that $K_0(S)\simeq \mathbf{Z}$.
Prove that any ring endomorphism $\varphi$ of $\mathbf{Z}\times\Gr$ in
$\Ho$ is of the form $\Psi^k$ for some $k\in\mathbf{Z}$. (Hint: $\varphi$
corresponds to a series $f\in\mathbf{Z}[[U]]$ which satisfies $f(0)=1$ and
$f(U)\cdot f(V)=f(U+V+UV)$. Then, ratiocinate in $\mathbf{Q}[[U]]$ to prove
that $f$ is of the form $(1+U)^\alpha$ for $\alpha\in\mathbf{Q}$.)
\end{exercise}

\subsection{The $\mathbf{P}^1$-spectrum $\BGL$}

Let $S$ be a regular scheme. We define a morphism $\sigma\colon
\mathbf{P}^1\wedge (\mathbf{Z}\times\Gr)\to\mathbf{Z}\times\Gr$ in $\Hopt$
(where $\infty$ is the base-point of $\mathbf{P}^1$) as the composition 
\[\xymatrix{\mathbf{P}^1\wedge (\mathbf{Z}\times\Gr)\ar[r]^-{u\wedge
\id}&
(\mathbf{Z}\times\Gr)\wedge (\mathbf{Z}\times\Gr)
\ar[r]^-{\mu} & \mathbf{Z}\times\Gr}\]
where $u\colon \mathbf{P}^1\to\mathbf{Z}\times\Gr$ corresponds to the class
$u=[\mathcal O(1)]-1\in \ker(\infty^\star\colon K_0(\mathbf{P}^1)\to
K_0(S))$ and $\mu$ is the pairing defined in
subsection~\ref{subsection-structures-higher-k-groups}. We denote
$\tilde{\sigma}\colon \mathbf{Z}\times\Gr\to
\R\SheafHom_\bullet(\mathbf{P}^1,\mathbf{Z}\times\Gr)$ the morphism in
$\Hopt$ corresponding to $\sigma$ by adjunction. It follows from the
projective bundle theorem that $\tilde{\sigma}$ is an isomorphism.

We can use this to define an object of the naive variant $\SHnaive$ (see
\cite[\S6]{riou-smf}) of the stable homotopy category $\SH$, \ie, an
($\Omega$)-$\mathbf{P}^1$-spectrum
up to homotopy. More precisely, an object of $\SHnaive$ consists in the
datum of a sequence $(\mathbf{E}_n)_{n\in\mathbf{N}}$ of objects of $\Hopt$
and of bonding morphisms $\sigma\colon \mathbf{P}^1\wedge \mathbf{E}_n\to
\mathbf{E}_{n+1}$ in $\Hopt$ which are supposed to be such that the
adjoint morphisms $\mathbf{E}_n\to
\R\SheafHom_\bullet(\mathbf{P}^1,\mathbf{E}_{n+1})$ are isomorphisms for
all $n\in\mathbf{N}$. The object $\BGLnaive\in\SHnaive$ is defined by
the fact that $(\BGLnaive)_n=\mathbf{Z}\times\Gr$ and that all bonding
morphisms identifies to the morphism $\sigma$ defined above.

We shall see that we may define an object $\BGL\in \SH$ up to a unique
isomorphism and that it lifts $\BGLnaive$. The obstruction we may encounter
to do this lies in the notion of stably phantom morphisms. More precisely,
if $\mathbf E$ and $\mathbf F$ are objects of $\SH$, represented by
$\Omega$-spectra, for any $i\in\mathbf{N}$, the sequence of groups
$(\Hom_{\Hopt}(\mathbf{E}_n,\R\Omega^i\mathbf{F}_n))_{n\in\mathbf{N}}$
is equipped with the structure of a projective system, and it follows from the
Milnor exact sequence that we have a short exact sequence (see
\cite[Lemme~6.5]{riou-smf}):

\begin{multline*}
0\to \R^1\lim_n \Hom_{\Hopt}(\mathbf{E}_n,\R\Omega\mathbf{F}_n)\to
\Hom_{\SH}(\mathbf E,\mathbf F)\\
\to \Hom_{\SHnaive}(\oub \mathbf{E},\oub \mathbf{F}) \to 0 \;\text{,}
\end{multline*}
where $\oub\colon \SH\to\SHnaive$ is the forgetful functor. The group on
the right identifies to $\lim_n\Hom_{\Hopt}(\mathbf{E}_n,\mathbf{F}_n)$ and
the group on the left is the subgroup of stably phantom morphisms
$\mathbf{E}\to\mathbf{F}$.

An object of $\SHnaive$ always lifts to an object of $\SH$, unique up to
isomorphism; however, this lifting is unique up to a \emph{unique}
isomorphism if and only if a given lifting has no nonzero stably phantom
endomorphisms \cite[Proposition~6.3]{riou-smf}. We will see that it is
the case for $\BGLnaive$ if $K_1(S)$ is finite (\eg, $S=\SpecZ$), which is
sufficient to construct a canonical $\BGL\in\SH$ for all regular schemes
$S$ as we may take the inverse image by $S\to \SpecZ$ (see
\cite[Proposition~4.4]{riou-smf}) of the unique one in $\SH[\SpecZ]$. 
This appeared in my thesis \cite{riou-these} and in
\cite{panin-pimenov-roendings} similar arguments reappeared.
This
being said, until the end of this subsection, we choose a lifting $\BGL$ of
$\BGLnaive$ in $\SH$.

\begin{remark}
In the study of projective systems
$(\Hom_{\Hopt}(\mathbf{E}_n,\R^1\Omega^i\mathbf{F}_n))_{n\in\mathbf{N}}$,
for some $i\in\mathbf{N}$, we may focus on the subsystem made of $H$-group
morphisms, which may be denoted
$(\Hom_{\Hopt}^+(\mathbf{E}_n,\R^1\Omega^i\mathbf{F}_n))_{n\in\mathbf{N}}$.
Indeed, the cokernel of this inclusion is a projective system with zero
transition maps, which implies that the inclusion induce isomorphisms on
$\lim$ and $\R^1\lim$.
\end{remark}

\begin{definition}
Let $A$ be an abelian group. We set $A^\Omega$ to be the following
projective system indexed by $\mathbf{N}$:
\[\dots \to
A[[U]]\vers{\Omega_{\mathbf{P}^1}}A[[U]]\vers{\Omega_{\mathbf{P}^1}}A[[U]]\vers{\Omega_{\mathbf{P}^1}}A[[U]]\;\text{,}\]
where the map $\Omega_{\mathbf{P}^1}\colon A[[U]]\to A[[U]]$ is defined by
$\Omega_{\mathbf{P}^1}(f)=(1+U)\frac{df}{dU}$.
\end{definition}

\begin{proposition}
Let $S$ be a regular scheme. The projective system
\[(\Hom_{\Hopt}^+((\BGL)_n,\R^i\Omega(\BGL)_n))_{n\in\mathbf{N}}\]
canonically identifies to $K_i(S)^\Omega$.
\end{proposition}

From proposition~\ref{proposition-splitting-principle}, we already know
that
$\Hom_{\Hopt}^+((\BGL)_n,\R\Omega^i(\BGL)_n)$ identifies degreewise
to the group
$K_i(S)[[U]]$. We let
$\omega\colon K_i(S)[[U]]\to K_i(S)[[U]]$ be the morphism corresponding to
the transition maps on the projective system
\[(\Hom_{\Hopt}^+((\BGL)_n,\R\Omega^i(\BGL)_n))_{n\in\mathbf{N}}\] under
this identification. We have to prove that $\omega=\Omega_{\mathbf{P}^1}$.

Let $\tau=\sum_{n\geq 0}a_nU^n\in K_i(S)[[U]]$. It corresponds to an
additive natural
transformation $(\tau_X\colon K_0(X)\to K_i(X))_{X\in\Sm}$ which is such
that $\tau_X([\mathcal L])=\sum_{n\geq 0} a_n([\mathcal L]-1)^n$ for all
line bundles $\mathcal L$.
The natural transformation $K_0(X)\to K_i(X)$ associated to $\omega(\tau)$
is characterised by the formula:
\[\omega(\tau)_X(x)\boxtimes v=\tau_{X\times \mathbf{P}^1}(x\boxtimes
v)\;\text{,}\]
where $v=[\mathcal O(1)]-1\in K_0(\mathbf{P}^1)$ and $\boxtimes$ is the
external product $K_\star(X)\times K_0(\mathbf{P}^1)\to K_\star(X\times
\mathbf{P}^1)$.
Assume that $x=[\mathcal L]$ is the class of a
line bundle $\mathcal L$ on a scheme $X\in \Sm$.
Then, $x\boxtimes v=
[\mathcal L\boxtimes \mathcal O(1)]-[\mathcal L\boxtimes \mathcal
O_{\mathbf{P}^1}]$. We may apply $\tau_{X\times \mathbf{P}^1}$ to this
difference; if we set $u=x-1$ and use that
$K_\star(X\times\mathbf{P}^1)\simeq K_\star(X)[v]/(v^2)$, we get:
\begin{eqnarray*}
\tau_{X\times\mathbf{P}^1}(x\boxtimes v)
&=&\sum_{n\geq 0}a_n\left[(1+u)(1+v)-1\right]^n-
\sum_{n\geq 0}a_n u^n\\
&=& \sum_{n\geq 0}a_n\left[(u+v(1+u))^n-u^n\right]\\
&=& \sum_{n\geq 1}na_n (1+u)u^{n-1}v\;\text{.}
\end{eqnarray*}
Then, $\omega(\tau)_X(x)=\sum_{n\geq 1}na_n(1+u)u^{n-1}$ which proves that
$\omega(\tau)=\sum_{n\geq
1}na_n(1+U)U^{n-1}=(1+U)\frac{d\tau}{dU}=\Omega_{\mathbf{P}^1}(\tau)$.

\begin{corollary}\label{corollary-milnor-sequence-bgl}
Let $S$ be a regular scheme. For all $i\in\mathbf{Z}$, we have a canonical
short exact sequence:
\[0\to \R^1\lim K_{i+1}(S)^\Omega\to \Hom_{\SH}(\BGL,\BGL[-i])\to
\lim K_i(S)^\Omega\to 0\;\text{.}\]
\end{corollary}

\begin{proposition}
Let $A$ be an abelian group. If $A$ is either finite or divisible, then
\[\R^1\lim A^\Omega=0\;\text{.}\]
\end{proposition}

If $A$ is divisible, the map
$\Omega_{\mathbf{P}^1}\colon A[[U]]\to A[[U]]$ is surjective. Hence, the
result is obvious in this case.

As a sequence of abelian groups $0\to A'\to A\to A''\to 0$ leads to a short
exact sequence of projective systems $0\to A'^\Omega\to A^\Omega\to
A''^\Omega\to 0$, a simple \emph{d\'evissage} reduces the case of
a finite abelian group $A$ to the special case
of $A=\mathbf{F}_p$ for a prime number $p$. Then, we are reduced to the
following lemma, which was suggested by Yves Andr\'e :

\begin{lemma}
Let $p$ be a prime number. We define $L_{\mathbf{F}_p}\subset
\mathbf{F}_p[[U]]$ as the subgroup of series $f=\sum_{n\geq 0}a_nU^n$
such that for all $k\in\mathbf{N}$, $\sum_{i=0}^{p-1}a_{kp+i}=0$.
Then,
\begin{itemize}
\item[(i)] The image of $\Omega_{\mathbf{P}^1}\colon \mathbf{F}_p[[U]]\to
\mathbf{F}_p[[U]]$ is $L_{\mathbf{F}_p}$;
\item[(ii)] If $f\in L_{\mathbf{F}_p}$, there exists a unique $g\in
L_{\mathbf{F}_p}$ such that $\Omega_{\mathbf{P}^1}(g)=f$;
\item[(iii)] The canonical map $\lim\mathbf{F}_p^\Omega\to
(\mathbf{F}_p^\Omega)_0$ induces a bijection $\lim\mathbf{F}_p^\Omega\simeq
L_{\mathbf{F}_p}$;
\item[(iv)] The projective system $L_{\mathbf{F}_p}$ satisfies
Mittag-Leffler condition. In particular, $\R^1\lim \mathbf{F}_p^\Omega=0$.
\end{itemize}
\end{lemma}

Let $f=\sum_{n\geq 0} a_nU^n$ and $g=\sum_{b\geq 0}b_nU^n$ be two elements
of $\mathbf{F}_p[[U]]$. The relation $\Omega_{\mathbf{P}^1}(g)=f$ is
equivalent to the equalities $nb_n +(n+1)b_{n+1}=a_n$ for all $n\geq 0$.
They can be restated as $nb_n=(-1)^{n-1}\sum_{k=0}^{n-1}a_k$ for all
$n\in\mathbf{N}$. It follows that $f$ is in the image of
$\Omega_{\mathbf{P}^1}$ if
and only if $\sum_{k=0}^{n-1}a_k=0$ whenever $p$ divides $k$, \ie, $f\in
L_{\mathbf{F}_p}$. Then, the relation $\Omega_{\mathbf{P}^1}(g)=f$
determines the coefficients $b_n$ for $p$ not dividing $n$ but says nothing
about the coefficients $b_{kp}$ for all $k\in\mathbf{N}$. There is a unique
possible choice for those so as to obtain $g\in L_{\mathbf{F}_p}$. We
have proved (i) and (ii). (iii) and (iv) immediately follow.

\begin{corollary}
\label{corollary-endomorphisms-bgl}
Let $S$ be a regular scheme. Let $i\in\mathbf{Z}$. If $K_{i+1}(S)$ is
finite or divisible, then
\[\Hom_{\SH}(\BGL,\BGL[-i])\simeq \lim_i K_i(S)^\Omega\;\text{.}\]
In particular, if $K_1(S)$ is finite (\eg, $S=\SpecZ$),
$\End_{\SH}(\BGL)\simeq \lim K_0(S)^\Omega$, $\BGL$ has no nonzero stably
phantom endomorphism in $\SH$ and thus $\BGLnaive\in\SHnaive$ lifts to an
object $\BGL\in\SH$ which is defined up to a unique isomorphism.
\end{corollary}

\begin{proposition}\label{proposition-omega-torsionfree}
Let $A$ be a torsionfree abelian group such that $\Hom(\mathbf{Q},A)=0$
(\eg, $A=\mathbf{Z}$).
Then, the map $\lim A^\Omega\to (A^\Omega)_0=A[[U]]$ is injective.
\end{proposition}

To prove this, it suffices to check that if $f\in A[[U]]$ is
such than $\Omega_{\mathbf{P}^1}(\Omega_{\mathbf{P}^1}(f))=0$, then
$\Omega_{\mathbf{P}^1}(f)=0$. Indeed, let $g=\Omega_{\mathbf{P}^1}(f)$. The
equality $\Omega_{\mathbf{P}^1}(g)=0$ implies that $g$ is constant, \ie,
$g\in A$. Then, we have $\frac{df}{dU}=\frac g{1+U}$ so that there exists
$h\in A$ such that $f=g\log (1+U)+h$. This series, which makes sense in
$(A\otimes_{\mathbf{Z}}\mathbf{Q})[[U]]$ does not lie in $A[[U]]$ unless
$g$ is in the image of a morphism $\mathbf{Q}\to A$. It follows that
$\Omega_{\mathbf{P}^1}(f)=g=0$.

\begin{remark}
Thanks to corollary~\ref{corollary-endomorphisms-bgl}, endomorphisms of
$\BGL$ in $\SH[\SpecZ]$ can be described as compatible families of
series in $\mathbf{Z}[[U]]$.
Proposition~\ref{proposition-omega-torsionfree} shows that this information
can be reduced to a single element in $\mathbf{Z}[[U]]$. However, I do not
know to which subgroup of $\mathbf{Z}[[U]]$ these endomorphisms correspond.
It obviously contains $1+U$ and $1/(1+U)$, which corresponds to 
the identity $\mathbf{\Psi}^1$ and the duality $\mathbf{\Psi}^{-1}$ (see
subsection~\ref{subsection-adams-operations}). According to
\cite{adams-clarke}, this group is strictly bigger and even uncountable!
\end{remark}

\subsection{Adams operations on $\BGLQ$}
\label{subsection-adams-operations}

The triangulated category $\SH$ may be localised so as to invert certain or
all primes. For instance, we may define $\SHQ$ as the full subcategory of
$\SH$ consisting of objects $A$ such that for any prime $p$, the
multiplication by $p$ on $A$ is an isomorphism. The left adjoint
$-_{\mathbf{Q}}\colon \SH\to \SHQ$ to this inclusion is called the
$\mathbf{Q}$-localisation functor. We let $\BGLQ$ be the image of $\BGL$ by
this functor. Then, for any finitely presented object $X$ of
$\SH$\;\footnote{An object $X$ in a triangulated category $\mathcal T$
where coproducts exist is finitely presented if the functor $\Hom_{\mathcal
T}(X,-)$ from $\mathcal T$ to the category of abelian groups commutes with
(infinite) coproducts. They constitute a triangulated subcategory
${\mathcal T}^{\pf}$ of $\mathcal T$.
In the case $\mathcal T=\SH$, $\SHpf$ is
the pseudo-abelian hull of the triangulated subcategory generated by objects
of the form $(\mathbf{P}^1)^{-n}\wedge U_+$ for $U\in\Sm$ (see
\cite[Proposition~1.2]{riou-cras-sw}).}, the canonical
map
$\Hom_{\SH}(X,\BGL)\otimes_{\mathbf{Z}}\mathbf{Q}\to\Hom_{\SH}(X,\BGLQ)$ is
a bijection and the methods used to obtain
corollaries~\ref{corollary-milnor-sequence-bgl} and
\ref{corollary-endomorphisms-bgl} give the following result:

\begin{corollary}
Let $S$ be a regular scheme. For all $i\in\mathbf{Z}$, we have a canonical
isomorphism:
\[\Hom_{\SHQ}(\BGLQ,\BGLQ[-i])\simeq \lim
(K_i(S)\otimes_{\mathbf{Z}}\mathbf{Q})^\Omega\;\text{.}\]
\end{corollary}

\begin{definition}\label{definition-psi-k-stable}
For all $k\in\mathbf{Z}-\{0\}$, we let
$\mathbf{\Psi}^k\in\End_{\SHQ}(\BGLQ)$ be the
endomorphism corresponding to the family
$(k^{-n}(1+U)^k)_{n\geq 0}\in \lim
\mathbf{Q}^\Omega$ (this family will also be denoted $\mathbf{\Psi}^k$).
\end{definition}

We obviously have the relations $\mathbf{\Psi}^k\circ \mathbf{\Psi}^{k'} =
\mathbf{\Psi}^{kk'}$. These Adams operations are constructed here with
$\mathbf{Q}$-coefficients, but it suffices to invert $k$ to define
$\mathbf{\Psi}^k$ (there might exist an obstruction to uniqueness in
$\R^1\lim K_1(S)[\frac 1 k]^\Omega$, in which case we may, as above,
construct it first on $\SpecZ$ and change the base).

To obtain a better understanding of the ring of endomorphisms of $\BGLQ$,
we focus on projective systems $A^\Omega$ in the case where $A$ is a
$\mathbf{Q}$-vector space:

\begin{definition}
Let $n\geq 0$. We define $p_n=\frac 1{n!} \log^n(1+U)\in\mathbf{Q}[[U]]$.
For any $\mathbf{Q}$-vector space $A$, we define an application
$\sigma\colon A^{\mathbf{N}}\to A[[U]]$ by the formula
\[\sigma((a_n)_{n\in\mathbf{N}})=\sum_{n=0}^\infty a_np_n\;\text{.}\]
\end{definition}

The infinite sum makes sense because the $U$-valuation of $p_n$ equals
$n$ and thus tends to $+\infty$.

\begin{lemma}\label{lemma-powers-of-log}
For any $\mathbf{Q}$-vector space $A$, the morphism $\sigma\colon
A^{\mathbf{N}}\to A[[U]]$ is an isomorphism of topological groups. If we
let $s\colon A^{\mathbf{N}}\to A^{\mathbf{N}}$ be the shift operator
$s((a_n)_{n\geq 0})=(a_{n+1})_{n\geq 0}$, we have the equality $\sigma
\circ s=\Omega_{\mathbf{P}^1}\circ \sigma$.
\end{lemma}

The topologies considered on $A^{\mathbf{N}}$ and $A[[U]]$ are the infinite
product topologies of the discrete topology on $A$. Then, the first
statement obviously follows from the fact that the $U$-valuation of $p_n$
is $n$. The second follows from the equalities
$\Omega_{\mathbf{P}^1}(p_n)=p_{n-1}$ for all $n\geq 1$ and
$\Omega_{\mathbf{P}^1}(p_0)=0$.

\begin{proposition}\label{proposition-definition-Sigma}
For any $\mathbf{Q}$-vector space, we may define $\Sigma\colon
A^{\mathbf{Z}}\to \lim A^\Omega$ by the formula
\[\Sigma((a_n)_{n\in\mathbf{Z}})=
(\sigma(a_n,a_{n+1},a_{n+2},\dots))_{n\geq 0}\;\text{,}\]
\ie, $\Sigma((a_n)_{n\in\mathbf{Z}})=\sum_{n\in\mathbf{Z}}a_n\pi_n$
where $\pi_n=(p_{n+k})_{k\geq 0}\in \lim \mathbf{Q}^\Omega$
(with $p_i$ set to zero for $i<0$).
\end{proposition}

It immediately follows from lemma~\ref{lemma-powers-of-log} which
identifies the projective system $A^\Omega$ to the projective system
\[\dots \vers{s} A^{\mathbf{N}}\vers{s} A^{\mathbf{N}}\vers{s}
A^{\mathbf{N}}\;\text{,}\]
whose projective limit is $A^{\mathbf{Z}}$.

\begin{remark}\label{remark-topology-on-end-bglq}
If $A=K_0(S)\otimes_{\mathbf{Z}}\mathbf{Q}$, a variant of 
proposition~\ref{proposition-splitting-principle} identifies $A[[U]]$ to
\[\End_{\Sm^\opp\Ab}(K_0(-)\otimes_{\mathbf{Z}}\mathbf{Q})\;\text{,}\]
so that the
composition law induces a law $\bigstar$ on $A[[U]]$ (see also
remark~\ref{remark-law-bigstar}). The operator $\Omega_{\mathbf{P}^1}$
defines an endomorphism of the ring $(A[[U]],+,\bigstar)$ so that $\lim
A^{\Omega}$ inherits a structure of a topological ring, which is, as a ring,
isomorphic to $\End_{\SHQ}(\BGLQ)$.
\end{remark}

\begin{proposition}\label{proposition-subring-end-bglq}
If $\mathbf{Q}^{\mathbf{N}}$ is endowed with its obvious ring structure and
$\mathbf{Q}[[U]]$ with the law $\bigstar$, then $\sigma\colon
\mathbf{Q}^{\mathbf{N}}\to \mathbf{Q}[[U]]$ is an isomorphism of
topological rings. The same conclusion applies to the isomorphism
$\Sigma\colon
\mathbf{Q}^\mathbf{Z}\isomto \lim \mathbf{Q}^\Omega$ whose target
identifies to a subring of $\End_{\SHQ}(\BGLQ)$ for any nonempty regular
scheme $S$.
\end{proposition}

We know that the $\mathbf{Q}$-vector space of $\mathbf{Q}[[U]]$ spanned by
elements $\Psi^k=(1+U)^k$, $k\geq 0$, is dense in $\mathbf{Q}[[U]]$. Hence,
it remains to prove the consistency of the formulas
$\Psi^{kk'}=\Psi^k\bigstar \Psi^{k'}$ with respect to the application of
$\sigma^{-1}\colon A[[U]]\isomto \mathbf{A}^{\mathbf{N}}$. This springs
from the following lemma:

\begin{lemma}
Let $k\in\mathbf{Z}-\{0\}$. Then,
\[(1+U)^k=\sigma((k^n)_{n\geq 0})\;\text{,}\quad
\mathbf{\Psi}^k=\Sigma((k^n)_{n\in\mathbf{Z}})\;\text{.}\]
\end{lemma}

Let $(\lambda_n)_{n\geq 0}=\sigma^{-1}(\Psi^k)$, where $\Psi^k$ is
identified to $(1+U)^k$. We know that
$\Omega_{\mathbf{P}^1}(\Psi^k)=k\Psi^k$. Then,
lemma~\ref{lemma-powers-of-log} implies that for all $n\geq 0$,
$\lambda_{n+1}=k\lambda_n$, so that $\lambda_n=k^n\lambda_0$. It remains to
compute $\lambda_0$. But, as it is the constant term of the series
$(1+U)^k$, we finally get $\lambda_0=1$.

\begin{definition}
For any $n\in\mathbf{Z}$, 
the element $\pi_n\in\lim\mathbf{Q}^\Omega$ was introduced in
proposition~\ref{proposition-definition-Sigma} and it is also the image
by $\Sigma\colon \mathbf{Q}^{\mathbf{Z}}\to \lim\mathbf{Q}^\Omega$ of the
characteristic function of $\{n\}\subset \mathbf{Z}$. Thanks to
proposition~\ref{proposition-subring-end-bglq}, for any regular scheme $S$,
$\pi_n$ identifies to an idempotent of $\End_{\SHQ}(\BGLQ)$. As $\SH$ has
infinite sums, it is pseudo-abelian (see
\cite[Proposition~II.1.2.9]{verdier}) and we may denote
$\BGLQ^{(n)}\subset \BGLQ$
the image of the projector $\pi_n$.
\end{definition}

\begin{theorem}
\label{theorem-diagonalisation-bglq}
Let $S$ be a regular scheme. The obvious morphism
\[\bigoplus_{n\in\mathbf{Z}}\BGLQ^{(n)}\to \BGLQ\]
is an isomorphism in $\SH$.
\end{theorem}

Let $n\geq 0$. We let $\chi_{[-n,n]}$ be the characteristic function of
$\{-n,\dots,n\}\subset \mathbf{Z}$. The corresponding element of $\lim
\mathbf{Q}^\Omega$ \emph{via} $\Sigma$ and the associated endomorphism of
$\BGLQ$ are also denoted $\chi_{[-n,n]}$. It is the sum of the orthogonal
idempotents $\pi_i$ for $-n\leq i\leq n$. Then, the image of
$\chi_{[-n,n]}$ identifies to $\bigoplus_{-n\leq k\leq n}\BGLQ^{(k)}$.

To prove that the morphism above is an isomorphism, it suffices
to prove that for any finitely presented object $X\in\SH$, the induced map
\[\Hom_{\SH}(X,\bigoplus_{n\in\mathbf{Z}}\BGLQ^{(n)})\to
\Hom_{\SH}(X,\BGLQ)\]
is a bijection. Due to previous observations, this map is injective and its
image is made of elements $x\in\Hom_{\SH}(X,\BGLQ)$ such that for a big
enough $n$, $\chi_{[-n,n]}(x)=x$. As the sequence 
$(\chi_{[-n,n]})_{n\in\mathbf{N}}$ of elements of $\mathbf{Q}^{\mathbf{Z}}$
tends pointwise to the constant function $\mathbf{1}$, the theorem shall be
a consequence of the following general lemma:

\begin{lemma}\label{lemma-comparison-topologies-end-bgl}
Let $S$ be a regular scheme. Let $(f_n)_{n\in\mathbf{N}}$ be a sequence of
elements in the group
$\lim (K_0(S)\otimes_{\mathbf{Z}}\mathbf{Q})^\Omega$ which
converges to an element $f$. Then, for any finitely presented object
$X$ in $\SH$ and $x\in\Hom_{\SH}(X,\BGLQ)$, there exists an integer $N$
such that for all $n\geq N$, $f_n(x)=f(x)$, where $f_n$ and $f$ are
identified to endomorphisms of $\BGLQ$.
\end{lemma}

Using the fact that the triangulated category $\SHpf$ identifies to the
pseudo-abelian hull of the category $\SWpf$ \cite[page
591]{voevodsky-icm}, we may assume that $X=(\mathbf{P}^1)^{\wedge -k}\wedge Y$
where $Y$ is a space of finite type (\eg, $S^i\wedge U_+$ where
$i\geq 0$ and $U\in\Sm$). Then, we are reduced to an unstable lemma:

\begin{lemma}
Let $S$ be a regular scheme. Let $(\tau_n)_{n\in\mathbf{N}}$ be a sequence
of elements in the group
$K_0(S)_{\mathbf{Q}}[[U]]$ which converges to an element
$\tau$. Then, for any space of finite type $X\in\Hopt$ and
$y\in\Hom_{\Hopt}(X,\mathbf{Z}\times\Gr)$, there exists $N\geq 0$ such that
for all $n\geq N$, $\tau_n(y)=\tau(y)$.
\end{lemma}

Variants of lemma~\ref{lemma-convergence-operations} and 
remark~\ref{remark-topology-operations} show that the lemma is true if $Y=U_+$
with $U\in\Sm$.  It holds more generally if $X$ is a pointed smooth
$S$-scheme,
for the obvious map
$\Hom_{\Hopt}(X,\mathbf{Z}\times\Gr)\to\Hom_{\Hopt}(X_+,\mathbf{Z}\times\Gr)$
is a split monomorphism, which, after tensoring with $\mathbf{Q}$, commutes
to $\tau$ and the $\tau_n$.

The general case
follows. As $X$ is of finite type, any
$x\in\Hom_{\Hopt}(Y,\mathbf{Z}\times\Gr)$ will factor through a disjoint
union of finite Grassmann varieties. Then, there exists a pointed smooth
$S$-scheme $U$, $u\in\Hom_{\Hopt}(U,\mathbf{Z}\times\Gr)$ and $f\colon
Y\to U$ in $\Hopt$ such that $y=f^\star(u)$. By the previous case, there
exists an integer $N$ such that $\tau_n(u)=\tau(u)$ for $n\geq N$. Hence,
$\tau_n(y)=\tau_n(f^\star u)=f^\star \tau_n(u)=f^\star \tau(u)=\tau(y)$ for
$n\geq N$.

\begin{remark}
One may find some inspiration from
lemma~\ref{lemma-comparison-topologies-end-bgl} so as to define a
topology on groups of morphisms $\Hom_{\mathcal T}(\mathbf E,\mathbf F)$
in a triangulated category $\mathcal T$ (where coproducts exist): the
weakest one such that for any morphism $x\colon X\to \mathbf E$ with
$X\in{\mathcal T}^\pf$, the
composition with $x$ induces a continuous map
$\Hom_{\mathcal T}(\mathbf E,\mathbf F)\to \Hom_{\mathcal T}(X,\mathbf F)$
where the target is endowed with the discrete topology. Then,
the lemma would say that in the case of $\End_{\SH}(\BGLQ)$,
this topology is the same as the one introduced in
remark~\ref{remark-topology-on-end-bglq}.
\end{remark}

\begin{proposition}
For any $n\in\mathbf{Z}$, the direct factor
$\BGLQ^{(n)}$ of $\BGLQ$ is preserved by $\mathbf{\Psi}^k$ for all
$k\in\mathbf{Z}-\{0\}$ and $\mathbf{\Psi}^k$ acts on it by multiplication
by $k^n$.
\end{proposition}

It follows from the following equalities in $\End_{\SH}(\BGLQ)$:
\[\mathbf{\Psi}^k\circ \pi_n=\pi_n\circ
\mathbf{\Psi}^k=k^n\pi_n\;\text{,}\] which can be proved
using their interpretations in the commutative subring
$\mathbf{Q}^\mathbf{Z}$ (see
proposition~\ref{proposition-subring-end-bglq}).

\begin{corollary}
For all $k\in\mathbf{Z}-\{0,\pm 1\}$ and $n\in\mathbf{Z}$, the endomorphism
$\mathbf{\Psi}^k-k^n\id$ of $\BGLQ$ has a kernel which is $\BGLQ^{(n)}$.
\end{corollary}

Using easy computations in $\mathbf{Q}^{\mathbf{Z}}$, we get the existence
of an automorphism $\phi_{n,k}$ of $\BGLQ$ such that $\phi_{n,k}\circ
(\mathbf{\Psi}^k-k^n\id)=\id-\pi_n$. Hence, the kernel of
$\mathbf{\Psi}^k-k^n\id$ is the same as the kernel of $\id-\pi_n$, which is
$\BGLQ^{(n)}$ by definition.

\bigskip

In other words, the decomposition of theorem~\ref{theorem-diagonalisation-bglq}
can be thought as a decomposition of $\BGLQ$ into a sum of eigenspaces
$\BGLQ^{(n)}$ for the Adams operations.

\begin{remark}
For any $a\in K_i(S)_{\mathbf{Q}}$, the constant family $a(1+U)$ belongs to
$\lim K_i(S)_{\mathbf{Q}}^\Omega$ (it can be interpreted as the natural
transformation $K_0(-)_\mathbf{Q}\to K_i(-)_\mathbf{Q}$ given by the
multiplication by $a$). It induces a morphism $\mu_a\colon
\BGLQ\to\BGLQ[-i]$. If $a\in K_i(S)^{(r)}$, one easily sees that
$\mu_a$ maps $\BGLQ^{(n)}$ to $\BGLQ^{(n+r)}[-i]$ for all $n\in\mathbf{Z}$.
Hence, we get a map
\[K_i(S)^{(r)}\to \Hom_{\SH}(\BGLQ^{(n)},\BGLQ^{(n+r)}[-i])\]
which is easily shown to be a bijection for all $n\in\mathbf{Z}$ and
$r\in\mathbf{Z}$.

It $S$ is a regular scheme of finite Krull dimension, the
$\gamma$-filtration on $K_i(X)$ has finitely many steps for all
$X\in\Sm$ (see \cite[\S2]{soule}); it can be used to prove that $\BGLQ$
is not only the direct sum of the $\BGLQ^{(n)}$ but also their infinite
product in $\SH$. This allows to give a description of morphisms
$\BGLQ\to\BGLQ[-i]$ as infinite matrices
$(a_{m,n})_{(m,n)\in\mathbf{Z}^2}$  where $a_{m,n}\in
K_i(S)^{(m-n)}$ corresponds to $\mu_{a_{m,n}}\colon
\BGLQ^{(n)}\to\BGLQ^{(m)}[-i]$.
\end{remark}

\begin{definition}\label{definition-h-beilinson}
Let $S$ be a regular scheme.
We set $\HBeilinson=\BGLQ^{(0)}\in\SHQ$.
\end{definition}

Using the periodicity isomorphism
$\SheafHom_\bullet(\mathbf{P}^1,\BGLQ)\simeq \BGLQ$, we get canonical
isomorphisms $\BGLQ^{(n)}\simeq \HBeilinson\wedge (\mathbf{P}^1)^{\wedge
n}$ for all $n\in\mathbf{Z}$.

\begin{remark}\label{remark-h-beilinson}
By its ``definition'' as an
eigenspace of Adams operations on the object $\BGLQ$ which
represents rationalized algebraic $K$-theory, this object $\HBeilinson$
represents motivic cohomology as it was first introduced by Beilinson (see
\cite{beilinson}).
\end{remark}

\section{Riemann-Roch theorems}
\label{section-riemann-roch}
\subsection{Adams-Riemann-Roch}

The Adams-Riemann-Roch theorem \cite[Theorem~7.6]{fulton-lang}
says that if $f\colon X\to S$ is a
projective morphism between regular schemes, then for all
$k\in\mathbf{Z}-\{0\}$ and $x\in K_0(X)\otimes_{\mathbf{Z}}\mathbf{Q}$:
\[\Psi^k(f_\star x)=f_\star (\Psi^k x \cdot (\theta_k
\Omega_f)^{-1})\;\text{,}\]
where $f_\star\colon K_0(X)\to K_0(S)$ is the direct image in $K$-theory
and $\theta_k\Omega_f$ is Bott's cannibalistic class associated to the
virtual cotangent bundle. It can be stated as a commutative square:
\[
\xymatrix@C=2cm{K_0(X) \ar[r]^{\Psi^k(-)\cdot (\theta_k
\Omega_f)^{-1}}\ar[d]^{f_\star} & K_0(X) \ar[d]^{f_\star} \\
K_0(S) \ar[r]^{\Psi^k} & K_0(S)}
\]

We shall obtain that for a projective and smooth morphism between regular
schemes, this diagram can be refined as a commutative diagram in $\SH$,
where $K_0(S)$ is replaced by $\BGLQS$ (we add the subscript
$S$ as a remainder of the base scheme) and $K_0(X)$ by $\R f_\star
\BGLQS[X]$ where $\R f_\star\colon \SH[X]\to \SH$ is the functor
constructed in \cite[Proposition~4.4]{riou-smf}. The proof will proceed by
showing that the diagram in $\SH$ commutes if and only if the relation
stated at the level of $K_0$ in the standard Adams-Riemann-Roch theorem is
true not only for $f\colon X\to S$ but for all morphisms
$f_T\colon X_T \to T$ deduced from $f$ by base change along smooth
morphisms $T\to S$.

One may expect that the homotopic version of Adams-Riemann-Roch we state
below (see theorem~\ref{theorem-adams-riemann-roch-sh}) has both sense and
truthfulness for more general projective morphisms between regular schemes.
However, the assumption that $f$ is projective and smooth shall be used
at several steps and thus should be considered as important in this
method.

\subsubsection{Pushforwards on $\BGL$}

\begin{proposition3}\label{proposition-f-star-stable}
Let $f\colon X\to S$ be a projective and smooth morphism between regular
schemes. There exists a morphism $\R f_\star
\BGLS[X]\vers{f_\star}
\BGLS$ in $\SH$ such that for any $n\in\mathbf{Z}$, $i\in\mathbf{N}$,
$T\in\Sm$, the map induced after applying the functor
$\Hom_{\SH}((\mathbf{P}^1)^{\wedge n}\wedge S^i\wedge T_+,-)$ identifies to
the usual pushforward in $K$-theory $f_\star\colon K_i(X_T)\to K_i(T)$ where
$X_T=X\times_S T$.
\end{proposition3}

\begin{lemma3}\label{lemma-f-star-unstable}
\sloppy Let
$f\colon X\to S$ be a projective and smooth morphism between regular
schemes. There exists a morphism $\R
f_\star(\mathbf{Z}\times\Gr_X)\vers{f_\star} \mathbf{Z}\times \Gr_S$ in
$\Hopt$ such that after the application of $\Hom_{\Hopt}(S^i\wedge T_+,-)$
for all $T\in \Sm$, we get the usual pushforward in $K$-theory
$f_\star\colon K_i(X_T)\to K_i(T)$ where $X_T=X\times_S T$.
\end{lemma3}

We have to use an homotopical description of these pushforwards in a way
which should be stricly functorial in $T\in \Sm$. We use Thomason's model
\cite[Lemma~3.5.3]{thomason-trobaugh}: for any regular scheme $X$, we
consider the complicial biWaldhausen category $\mathcal C(X)$ of perfect
bounded above complexes of flat $\mathcal O_X$-modules\;\footnote{Note that
we have to fix suitable cardinality bounds so as to get (essentially) small
categories.}. For any (regular) base scheme $S$, it is easy to turn this
construction into a presheaf $\mathcal C_S$ of complicial biWaldhausen
categories over $\Sm$, with $\mathcal C_S(X)$ equivalent to $\mathcal C(X)$
for all $X\in \Sm$. Then, the associated presheaf of $K$-theory spaces
$\K\mathcal C_S$ is a model of algebraic $K$-theory (\ie, it is
canonically isomorphic to $\mathbf{Z}\times\Gr$ in $\Hopt$, see
Proposition~\ref{proposition-genuine-models}).
At this stage, it is obvious that $\K\mathcal C_X$ is acyclic for the
functor $f_\star$, \ie, we have a canonical isomorphism
$\R f_\star \K\mathcal C_X\simeq f_\star \K\mathcal C_X$ in $\Hopt$.

We shall construct the expected morphism $f_\star\colon \R f_\star
(\mathbf{Z}\times\Gr_X)\to \mathbf{Z}\times\Gr_S$ as a morphism $f_\star
\K\mathcal C_X\to\K\mathcal C_S$. The details follow. We choose a finite
open cover $\mathcal U=\{U_1,\dots,U_n\}$ of $X$ such that all the induced
morphisms $f_i\colon U_i\to S$ are affine (as we assumed $S$ separated, any
affine open cover of $X$ has this property). For any nonempty subset $I$ of
$\{1,\dots,n\}$, we set $U_I=\cap_{i\in I}U_i$ and denote
$f_I\colon U_I\to S$ the restriction of $f$ to these subschemes.

For any $T\in \Sm$, we consider the base change $f_T\colon X_T\to T$ of $f$
along $T\to S$ and introduce the morphisms $f_{I,T}\colon U_I\times_S T\to
T$ deduced from $f_I$ for all nonempty subsets $I$ of $\{1,\dots,n\}$.
These morphisms $f_{I,T}$ are affine and flat. For any $\mathcal M\in
\mathcal C(X_T)$, we define $(f_{\bullet,T})_\star \mathcal M$
as the total complex of
the \v{C}ech type bicomplex:
\[\dots \to 0\to \oplus_{1\leq i\leq n} (f_{i,T})_\star \mathcal M\to
\oplus_{1\leq i<j\leq n} (f_{i,j,T})_\star \mathcal M\to 
\dots\;\text{,}\]
where the first \emph{a priori} non trivial object lies in cohomological
degree $0$. As $f$ is flat,
the object $(f_{\bullet,T})_\star \mathcal M$ is a bounded
complex of flat $\mathcal O_T$-modules
and from standard results in coherent cohomology (see
\cite[Th\'eor\`eme 3.2.1]{EGA3Vol1}),
$(f_{\bullet,T})_\star\mathcal M$ represents $\R {f_T}_\star \mathcal M$ in
the derived category $D(T,\mathcal O_T)$ and is perfect. Hence, we have
defined a functor $(f_{\bullet,T})_\star
\colon \mathcal C(X_T)\to \mathcal C(T)$ for any $T\in\Sm$. This
construction commutes up to canonical isomorphisms with the inverse image
functors (\ie, the presheaf structure on $\mathcal C_S$)
associated to morphisms $T'\to T$ in $\Sm$. It is an easy game to
modify the definitions so as to get strict compatibilities. Finally, we may
apply the $K$-theory functor to obtain the expected morphism $f_\star
\K\mathcal C_X\vers{f_\star} \K\mathcal C_S$ of presheaves of pointed sets
on $\Sm$.

\medskip

The compatibility between pushforwards and external products implies that
we may use the morphism from lemma~\ref{lemma-f-star-unstable} to
define a morphism $\R f_\star \BGLS[X]\to \BGLS$ up to stably phantom maps
(\ie, in $\SHnaive$). In the statement of
proposition~\ref{proposition-f-star-stable}, there is no uniqueness claim.
However, we shall see in the sequel that it will be the case after
tensoring with $\mathbf{Q}$.

\subsubsection{Statement of the theorem}

\begin{theorem3}\label{theorem-adams-riemann-roch-sh}
Let $f\colon X\to S$ be a projective and smooth morphism between regular
schemes. Then, the following diagram in $\SH$ commutes:
\[
\xymatrix@C=3cm{
\R f_\star \BGLQS[X]\ar[d]^{f_\star}
\ar[r]^{\R f_\star(\theta_k(\Omega_f)^{-1}\cdot \mathbf{\Psi}^k(-))} &
\R f_\star \BGLQS[X]
\ar[d]^{f_\star} \\
\BGLQS \ar[r]^{\mathbf{\Psi}^k}\ar[r] & \BGLQS
}
\]
where both vertical maps are the pushforward morphism constructed in
proposition~\ref{proposition-f-star-stable} (tensored with $\mathbf{Q}$),
the lower map is $\mathbf{\Psi}^k\in \End_{\SH}(\BGLQS)$
(see definition~\ref{definition-psi-k-stable}) and the upper map is
obtained by applying $\R f_\star$ to the endomorphism of $\BGLQS[X]$
corresponding to $\mathbf{\Psi}^k$ multiplied by the inverse of Bott's
cannibalistic class\;\footnote{It makes sense as previous results show that
$\End_{\SH[X]}(\BGLQS[X])$ is a module over
$K_0(X)\otimes_{\mathbf{Z}}\mathbf{Q}$.}.
\end{theorem3}

\begin{corollary3}\label{corollary-adams-riemann-roch-k-i}
Let $f\colon X\to S$ be a projective and smooth morphism between regular
schemes. Then, the following diagram commutes for any $i\in\mathbf{Z}$:
\[
\xymatrix@C=2cm{K_i(X) \ar[r]^{\Psi^k(-)\cdot (\theta_k
\Omega_f)^{-1}}\ar[d]^{f_\star} & K_i(X) \ar[d]^{f_\star} \\
K_i(S) \ar[r]^{\Psi^k} & K_i(S)}
\]
\end{corollary3}

Corollary~\ref{corollary-adams-riemann-roch-k-i} is deduced from the
statement of theorem~\ref{theorem-adams-riemann-roch-sh} by applying
functors $\Hom_{\SH}(S^i\wedge T_+,-)$. Conversely, I claim that two
morphisms $\R f_\star \BGLQS[X]\to \BGLQS$ in $\SH$ are equal as soon as
they induce equal maps after the application of functors
$\Hom_{\SH}((\mathbf{P}^1)^{\wedge -n}\wedge T_+,-)$ for all
$n\in\mathbf{N}$ and $T\in \Sm$. This will be the goal of
theorem~\ref{theorem-morphisms-f-star-bgl-to-bgl} in the paragraph which
follows. Then, theorem~\ref{theorem-adams-riemann-roch-sh} shall follow
from the classical Adams-Riemann-Roch theorem (\ie, the case $i=0$
in corollary~\ref{corollary-adams-riemann-roch-k-i}).

\subsubsection{Morphisms $\R f_\star \BGLQS[X]\to \BGLQS$}
\label{subsubsection-morphisms-f-star-bgl-to-bgl}

\begin{definition3}
For all $(i,j)\in\mathbf{Z}^2$, we define a functor $\pi_{i,j}\colon
\SH\to\Sm^\opp\Ab$ by 
\[(\pi_{i,j}\mathbf{E})(U)=\Hom_{\SH}((\mathbf{P}^1)^{\wedge
j}\wedge S^{i-2j}\wedge U_+,\mathbf{E})\;\text{;}\]
they are the functors ``presheaves of stable homotopy groups''.
\end{definition3}

\begin{theorem3}\label{theorem-morphisms-f-star-bgl-to-bgl}
Let $f\colon X\to S$ be a projective and smooth morphism between regular
schemes. Let $\tau\colon \R f_\star \BGLQS[X]\to \BGLQS$ be a morphism in
$\SH$ such that for
all $n\in\mathbf{Z}$, $\pi_{2n,n}(\tau)=0$\;\footnote{One may notice that
$\pi_{2n,n}(\tau)=0$ implies $\pi_{2(n+1),n+1}(\tau)=0$.}. Then, $\tau=0$.
\end{theorem3}

We use the theory of stable homotopic functors (see
\cite{ayoub-asterisque-i} and also \cite[Remarque~4.6]{riou-smf}). Thus, we
have a direct image functor with proper support $\R f_!\colon \SH[X]\to
\SH$ which has a right adjoint $f^!$. As $f$ is projective, we have a
canonical isomorphism $\R f_!\isomto \R f_\star$. Then, by adjunction, the
morphism $\tau\colon \R f_\star \BGLQS[X]\to \BGLQS$ corresponds to a
morphism $\tilde{\tau}\colon \BGLQS[X]\to f^!\BGLQS$.

\begin{lemma3}\label{lemma-thom-isomorphism-bgl} We let $f\colon X\to S$
be a projective and smooth morphism between regular schemes.

\begin{itemize}
\item[(i)]
There exists a canonical isomorphism $f^!\BGLQS\simeq \BGLQS[X]$ in
$\SH[X]$.
\item[(ii)] For any vector bundle $\mathcal E$ over $X$, we have a canonical
isomorphism $\BGLQS[X]\wedge \Th \mathcal E\simeq \BGLQS[X]$ in $\SH[X]$.
\end{itemize}
\end{lemma3}

By definition of $f^!$, for any $\mathbf{E}\in\SH$,
we have an isomorphism $f^!\mathbf{E}\simeq \L f^\star
\mathbf{E}\wedge \Th T_f$ where $T_f$ is the relative tangent bundle of $f$
and $\Th T_f$ its Thom space. As $\L f^\star \BGLQS$ identifies to
$\BGLQS[X]$, (i) will follow from (ii).

To prove (ii), we consider the isomorphism $\Th \mathcal E\simeq
\mathbf{P}(\mathcal E\oplus
\mathcal O_X)/\mathbf{P}(\mathcal E)$\;\footnote{We, reluctantly, do not follow
Grothendieck's convention. Here, $\mathbf{P}(\mathcal E)$ is the projectivisation of
the symmetric algebra of the dual of $\mathcal E$.} and the class $\xi$ of the
fundamental sheaf $\mathcal O(1)$ in $K_0(\mathbf{P}(\mathcal E\oplus \mathcal
O_X))$. We may set $v=\xi^r-[\wedge^1 \mathcal E]\xi^{r-1}+[\wedge^2
\mathcal E]\xi^{r-2}+\dots+(-1)^r[\wedge^r \mathcal E]\in
K_0(\mathbf{P}(\mathcal E\oplus \mathcal
O_X))$ where $r$ is the rank of $\mathcal E$.
The class $v$ vanishes when restricted
to $\mathbf{P}(\mathcal E)$. Hence, $v$ actually defines an element in
$\tilde{K}_0(\Th \mathcal E)$. In this paragraph, $\tilde{K}_0(Y)$ is the reduced
$K$-theory of a pointed space $Y$, \ie,
$\Hom_{\Hopt}(Y,\mathbf{Z}\times\Gr)$, which identifies to the kernel of
the map $K_0(Y)\to K_0(S)$ given by the base-point. Even if we use the same
notation, it should not be confused with the kernel of $\rk\colon K_0(X)\to
\mathbf{Z}^{\pi_0(X)}$, which makes sense for $X\in\Sm$.

Using the multiplicative structure on $\mathbf{Z}\times \Gr$, we may
consider the external product with $v$ in $\Hopt$:
\[\mathbf{Z}\times \Gr\to \R\SheafHom_\bullet(\Th
\mathcal E,\mathbf{Z}\times\Gr)\;\text{,}\]
which is seen to be an isomorphism thanks to computations using the 
projective bundle formula. Using this morphism termwise, we get the
expected isomorphism \[\BGLS[X]\isomto \R\SheafHom_\bullet(\Th
\mathcal E,\BGLS[X])\]
in $\SHnaive[X]$\;\footnote{This construction may also be deduced from a
more universal pairing $(\mathbf{Z}\times\Gr)\wedge \BGL\to \BGL$
which should be constructed first in $\SH[\SpecZ]$. However, 
when one want to tackle the trouble of stably phantoms
morphisms, one has to use different arguments than those appearing in this
article. To do this, we can use \cite[Lemma~A.6]{riou-these} which we used
there to obtain another proof of the construction of $\BGL$, see
corollary~\ref{corollary-endomorphisms-bgl}. This method can be continued in
order to obtain an associative and commutative
pairing $\BGL\wedge \BGL\to \BGL$ in $\SH$ (see
\cite[Theorem~2.2.1]{panin-pimenov-roendings}).}.
As $\Th \mathcal E$ is invertible for the $\wedge$-product on
$\SH[X]$ (see proposition~\ref{proposition-thom-spectrum}),
property~(ii) follows.

\begin{lemma3}\label{lemma-vanishing-on-negative-thom-spaces}
Let $\psi\colon \mathbf{E}\to \mathbf{F}$ be a morphism in $\SH$. We assume
that $\mathbf{F}$ is such that for any $U\in\Sm$, vector bundle $\mathcal
E$ on $U$ and $n\in\mathbf{Z}$, the canonical map
$\tilde{\mathbf{F}}^{2n,n}(\Th_U \mathcal E)\to
\mathbf{F}^{2n,n}(\mathbf{P}(\mathcal E\oplus \mathcal O_U))$ is
injective\;\footnote{We use standard implicit convention. More precisely,
this map is the result of the application of the functor
$\Hom_{\SH}((\mathbf{P}^1)^{\wedge n}\wedge -,\mathbf{F})$ to the canonical
morphism $\mathbf{P}(\mathcal E\oplus \mathcal O_U)_+\to \Th_U \mathcal E$
in $\Hopt$.}. We also assume that $\pi_{2n,n}(\psi)=0$ for all
$n\in\mathbf{Z}$. Then, for any vector bundle $\mathcal E$ on $U\in \Sm$
and any $n\in\mathbf{Z}$, the map $\Hom_{\SH[X]}((\mathbf{P}^1)^{\wedge
n}\wedge \Th_U(-\mathcal E),\psi)$ vanishes (see
definition~\ref{definition-thom-spectrum}).
\end{lemma3}

Using Jouanolou's trick, we may assume that $U$ is affine. Then, the virtual
bundle $-\mathcal E$ identifies to a difference $\mathcal F-\mathcal O_U^k$
where $\mathcal F$ is a genuine vector bundle and $k\in\mathbf{N}$. Then,
we want to prove that for any $n\in\mathbf{Z}$, the morphism
$\Hom_{\SH}((\mathbf{P}^1)^{\wedge n}\wedge \Th_U \mathcal F,\psi)$
vanishes. As, $\Th_U \mathcal F = \mathbf{P}(\mathcal F\oplus \mathcal
O_U)/\mathbf{P}(\mathcal F)$, the result follows from the second assumption
for $U=\mathbf{P}(\mathcal F\oplus \mathcal O_U)$ and the injectivity
stated in the first assumption.

\bigskip

Now, we shall prove theorem~\ref{theorem-morphisms-f-star-bgl-to-bgl}. We
may apply lemma~\ref{lemma-vanishing-on-negative-thom-spaces} to
$\tau\colon \R f_\star \BGLQS[X]\to \BGLQS$. Then, as 
$f_!\mathbf{E}\simeq \L f_\sharp(\mathbf{E}\wedge \Th_X(-Tf))$ for any
$\mathbf{E}\in\SH[X]$, we obtain the vanishing of the maps
$\Hom_{\SH}(f_!((\mathbf{P}^1)^{\wedge n}
\wedge U_+),\tau)$ for all $n\in\mathbf{Z}$
and $U\in\Sm[X]$. By adjunction, it implies that the maps
$\Hom_{\SH[X]}((\mathbf{P}^1)^{\wedge n}\wedge U_+,\tilde{\tau})$ vanish.
As we know that $\tilde{\tau}$ can be identified to an endomorphism of
$\BGLQS[X]$ (see lemma~\ref{lemma-thom-isomorphism-bgl}), we can use the
results of section~\ref{section-additive-and-stable-results} to assert
that $\tilde{\tau}=0$. Finally, by adjunction, $\tau=0$.

\subsection{Motivic Eilenberg-Mac Lane spectra}
\subsubsection{Morphisms $\mathbf{Z}\times\Gr\to K(\mathbf{Z}(n),2n)$}

\begin{definition3}
Let $k$ be a perfect field. For any $n\geq 0$, we denote
$K(\mathbf{Z}(n),2n)$ the motivic Eilenberg-Mac Lane space defined in
\cite[\S6.1]{voevodsky-icm}. For $i\geq 0$. we let $K(\mathbf{Z}(n),2n-i)$
be its $i$th loop space.
\end{definition3}

By definition, for any $n\geq 0$ and $i\geq 0$, the group
$\Hom_{\Ho[k]}(X,K(\mathbf{Z}(n),2n-i))$ identifies to the motivic cohomology
group $H^{2n-i}(X,\mathbf{Z}(n))$. The comparison with (higher) Chow groups
\cite{voevodsky-chow} implies that for any $n\geq 0$, there is a canonical
isomorphism $\pi_0K(\mathbf{Z}(n),2n)\simeq CH^n(-)$ in $\Sm[k]^\opp\Ab$ where $CH^n(-)$ is the presheaf $X\longmapsto CH^n(X)$.

\begin{theorem3}\label{theorem-morphisms-k-to-ch}
Let $k$ be a perfect field. Let $n\geq 0$. Then, the functor $\pi_0$
induces a bijection:
\[\Hom_{\Ho[k]}(\mathbf{Z}\times\Gr,K(\mathbf{Z}(n),2n))\isomto
\Hom_{\Sm[k]^\opp\Sets}(K_0(-),CH^n(-))\;\text{.}\]
Moreover, the graded algebra
$(\Hom_{\Sm[k]^\opp\Sets}(\tilde{K}_0(-),CH^n(-))_{n\in\mathbf{N}}$
identifies to the polynomial algebra $\mathbf{Z}[c_1,c_2,\dots]$ where
$c_i$ lies in degree $i$ and corresponds to the $i$th Chern class
$c_i\colon \tilde{K}_0(-)\to CH^i(-)$.
\end{theorem3}

The first statement follows from the fact that whenever $d\leq d'$ and
$r\leq r'$, the inclusion $\Gr_{d,r}\subset \Gr_{d',r'}$ induces a split
monomorphisms $M(\Gr_{d,r})\subset M(\Gr_{d',r'})$ on motives. This fact
follows from the cellularity of $\Gr_{d,r}$, $\Gr_{d',r'}$ and
$\Gr_{d',r'}-\Gr_{d,r}$ (see \cite[\S3]{kahn-cellular}
for a similar statement). Then,
any object representing a cohomology which factors through the category of
motives will satisfy property~(K) with any number of operands (see
definition~\ref{definition-property-k}) and we may use
theorem~\ref{theorem-property-k-implies-computation}.

The second part arises from the computation of Chow groups of Grassmann
varieties $\Gr_{d,r}$ for $d,r\geq 0$ (see \cite{grothendieck-chevalley})
and the passage to the limit $r\to \infty$ and $d\to \infty$ as it was done
for the algebraic $K$-theory.

\subsubsection{Additive morphisms}
\label{subsubsection-additive-morphisms-chow}

The proof of theorem~\ref{theorem-morphisms-k-to-ch} applies not only to
natural transformations $K_0(-)\to CH^n(-)$ but also to natural
transformations involving several operands, \eg, $K_0(-)\times K_0(-)\to
CH^n(-)$. Hence, $H$-group morphisms $\mathbf{Z}\times\Gr\to
K(\mathbf{Z}(n),2n)$ correspond to morphisms $K_0(-)\to CH^n(-)$ in
$\Sm[k]^\opp\Ab$ (see
proposition~\ref{proposition-lifting-algebraic-structures}). The group of
these morphisms is described in the following proposition:

\begin{proposition3}
Let $k$ be a perfect field. For any $n\geq 0$, the map given by the
evaluation at $[\mathcal O(1)]$ in $K_0(\mathbf{P}^n)$ induces an
isomorphism:
\[\Hom_{\Sm[k]^\opp\Ab}(K_0(-),CH^n(-))\isomto \lim_{r\in\mathbf{N}}
CH^n(\mathbf{P}^r)\simeq CH^n(\mathbf{P}^n)\simeq \mathbf{Z}\;\text{.}\]
We denote $\chi_n\colon K_0(-)\to CH^n(-)$ the canonical generator given
by this isomorphism. It is characterised by the fact that
$\chi_n([\mathcal L])=[D]^n$
anytime $\mathcal L$ is a line bundle on $X\in\Sm$ and $D$ is the divisor of a
rational section of $\mathcal L$.
\end{proposition3}

The proof of the injectivity of the map 
\[\Hom_{\Sm[k]^\opp\Ab}(K_0(-),CH^n(-))\isomto \lim_{r\in\mathbf{N}}
CH^n(\mathbf{P}^r)\]
is similar to that of proposition~\ref{proposition-splitting-principle}.
The group $\lim_{r\in\mathbf{N}} CH^n(\mathbf{P}^r)$ is easily identified
to the group $\mathbf{Z}$, generated by the compatible family made of $n$th
powers of classes in hyperplanes in $\mathbf{P}^r$ for all
$r\in\mathbf{N}$. For the surjectivity, \ie, the existence of $\chi_n$, we
shall use the following lemma, which is a consequence of the theory of
symmetric polynomials (hint: use \cite[VI~4.3-4.4]{SGA6}):

\begin{lemma3}
Let $n\geq 1$. There exists a unique functorial homomorphism 
\[\chi_n\colon (1+A[[t]]^+,\times)\to (A,+)\] for all commutative
rings $A$ such that for any $x\in A$,
\[\chi_n(1+xt)=x^n\;\text{,}\]
and $\chi_n$ vanishes on the subgroup $1+t^{n+1}A[[t]]$.
\end{lemma3}

Note that by looking at the universal situation, we know that
$\chi_n(\sum_{i\geq 0}a_it^i)$ is given by a polynomial in $a_1,\dots,a_n$
and it is homogeneous of total degree $n$ if we set $\deg a_i=i$.

For any $X\in\Sm[k]$, $u\in K_0(X)$, we consider the Chern polynomial
$c_t(u)\in CH^\star(X)[[t]]$ and apply the construction of the lemma to
this series : $\chi_n(c_t(u))\in CH^n(X)$. This constructs a natural
transformation $K_0(-)\to CH^n(-)$ to which we give the same name $\chi_n$.
This finishes the proof of the proposition in the case $n\geq 1$; the
remaining case $n=0$ is trivial.

\begin{remark3}
As we have seen it,
the natural transformation $\chi_n\colon K_0(-)\to CH^n(-)$ is given by
a polynomial involving Chern classes. It can be computed inductively using
Newton relations:
\[\chi_k-c_1\chi_{k-1}+\dots+(-1)^{k-1}c_{k-1}\chi_1+(-1)^kkc_k=0\;\text{.}\]
For instance, $\chi_1=c_1$, $\chi_2=c_1^2-2c_2$,
$\chi_3=c_1^3-3c_1c_2+3c_3$.
\end{remark3}

The following similar result gives a computation of the group of $H$-group
morphisms $\mathbf{Z}\times\Gr\to K(\mathbf{Z}(n),2n)$ in $\Hopt[k]$.

\begin{corollary3}\label{corollary-morphisms-z-gr-k-z-n-2n-i}
Let $k$ be a perfect field, $n\geq 0$, $i\geq 0$. For any $0\leq j\leq
\min(i,n)$ and $x\in H^{2j-i}(k,\mathbf{Z}(j))$, we define a natural
transformation
$x\cdot \chi_{n-j}\colon K_0(-)\to H^{2n-i}(-,\mathbf{Z}(n))$ of
presheaves of abelian groups on $\Sm[k]$, obtained as the composition of
$\chi_{n-j}$ and the multiplication by $x$ on motivic cohomology. Then, the
group of natural transformations $K_0(-)\to H^{2n-i}(-,\mathbf{Z}(n))$
identifies to the direct sum of the groups
$H^{2j-i}(k,\mathbf{Z}(j))$ for $0\leq j\leq \min(i,n)$, as follows:
\[\Hom_{\Sm[k]^\opp\Ab}(K_0(-),H^{2n-i}(-,\mathbf{Z}(n)))\simeq 
\bigoplus_{j=0}^{\min(i,n)}
H^{2j-i}(k,\mathbf{Z}(j))\cdot\chi_{n-j}\;\text{.}\]
\end{corollary3}

\subsubsection{Stable morphisms}

The motivic Eilenberg-Mac Lane spectrum $\HMot$ is obtained from the
sequence of objects $K(\mathbf{Z}(n),2n)$ (see \cite[\S
6.1]{voevodsky-icm}). We may describe its image in $\SHnaive[k]$ by saying
that the different Eilenberg-Mac Lane spaces are related by the canonical
isomorphism $K(\mathbf{Z}(n),2n)\simeq
\R\SheafHom_\bullet(\mathbf{P}^1,K(\mathbf{Z}(n+1),2n+2)$ induced by the
external product with the class of the $1$-codimensional cycle $[\infty]$ in
$CH^1(\mathbf{P}^1)$. This construction generalises to give a
$\mathbf{P}^1$-spectrum $\HMot[A]$ for any coefficient abelian group $A$.

In order to study morphisms $\BGL\to \HMot{}[-i]$ for $i\geq 0$, we use the
following definition.

\begin{definition3}
Let $n\geq 1$ and $i\geq 0$. Let $\tau\colon K_0(-)\to
H^{2n-i}(-,\mathbf{Z}(n))$ be an additive natural transformation,
\ie, a morphism in $\Sm[k]^\opp\Ab$. We define a natural
transformation $\Omega_{\mathbf{P}^1}(\tau)\colon K_0(-)\to
H^{2n-2-i}(-,\mathbf{Z}(n-1))$ which shall be characterised by the
commutativity of the following diagram for all $X\in\Sm[k]$:
\[
\xymatrix@C=2cm{K_0(X)\ar[r]^{u\boxtimes -} \ar[d]^{\Omega_{\mathbf{P}^1}(\tau)}
  & K_0(\mathbf{P}^1\times
X)\ar[d]^\tau \\
H^{2n-2-i}(X,\mathbf{Z}(n-1)) \ar[r]^{[\infty]\boxtimes -} & 
H^{2n-i}(\mathbf{P}^1\times X,\mathbf{Z}(n))
}
\]
where $u=[\mathcal O(1)]-1\in K_0(\mathbf{P}^1)$ and $[\infty]$ is the
class of a rational point in
$CH^1(\mathbf{P}^1)=H^2(\mathbf{P}^1,\mathbf{Z}(1))$.
\end{definition3}

\begin{lemma3}
\label{lemma-omega-chi-n}
Let $k$ be a perfect field. For any $n\geq 1$, we have
$\Omega_{\mathbf{P}^1}(\chi_n)=n\chi_{n-1}$.
\end{lemma3}

By the splitting principle, it suffices to check that
$\Omega_{\mathbf{P}^1}(\chi_n)$ and $n\chi_{n-1}$ coincide on elements of
the form $[\mathcal L]\in K_0(X)$ where $\mathcal L$ is a line bundle on
some $X\in\Sm[k]$. Let $D$ be the divisor of a rational section of
$\mathcal L$.
Considering $CH^\star(X\times\mathbf{P}^1)$ both as an algebra over
$CH^\star(X)$ and $CH^\star(\mathbf{P}^1)$, we get:
\begin{eqnarray*}
[\infty]\boxtimes \Omega_{\mathbf{P}^1}(\chi_n)([\mathcal L])&=&
\chi_n(u\boxtimes
[\mathcal L]) \\
&=&\chi_n([\mathcal O(1)\boxtimes \mathcal L])-\chi_n(\mathcal
O_{\mathbf{P}^1}\boxtimes \mathcal
L)\\
&=&([\infty]+[D])^n-[D]^n=n[\infty][D]^{n-1}\\
&=&[\infty]\boxtimes
(n\chi_{n-1}([\mathcal L]))\;\text{,}
\end{eqnarray*}
which proves the expected result: $\Omega_{\mathbf{P}^1}(\chi_n)([\mathcal
L])=n\chi_{n-1}([\mathcal L])$.

\bigskip

This lemma leads to a description of the projective system
\[(\Hom^+_{\Hopt}(\mathbf{Z}\times\Gr,K(\mathbf{Z}(n),2n)))_{n\in\mathbf{N}}\]
deduced from the bonding morphisms on $\BGL$ and $\HMot$; it identifies to
a projective system which we shall denote $\mathbf{Z}!$:
\[\dots\to \mathbf{Z}\vers 5 \mathbf{Z}\vers 4 \mathbf{Z}\vers 3  \mathbf{Z}\vers 2  \mathbf{Z}\vers 1
\mathbf{Z}\;\text{.}\]
We generalise this definition:

\begin{definition3}
Let $A$ be an abelian group. We define a projective system $A!$ of abelian
groups indexed by $\mathbf{N}$ by saying that in degree $n\in\mathbf{N}$,
$(A!)_n=A$ and the transition map $(A!)_n\to (A!)_{n-1}$ is the
multiplication by $n$ on $A$.
\end{definition3}

\begin{definition3}
If $X_\bullet=(\dots \to X_n\vers{f_{n-1}} X_{n-1}\to \dots \to X_1\vers{f_0}
X_0)$ is a projective system of abelian groups indexed by $\mathbf{N}$, we
define a new projective system $sX_\bullet =
(\dots \to X_n\vers{f_{n-1}} X_{n-1}\to \dots \to X_1\vers{f_0} X_0\to 0)$.
\end{definition3}

\begin{proposition3}
Let $A$ be an abelian group. We let $\HMot[A]$ be the motivic Eilenberg-Mac
Lane spectrum with coefficients in $A$. Then, for any $i\in\mathbf{Z}$,
the projective system
\[(\Hom^+_{\Hopt[k]}(\mathbf{Z}\times\Gr, K(A(n),2n-i))_{n\in\mathbf{N}}\]
associated to the $\mathbf{P}^1$-spectra $\BGL$ and $\HMot[A][-i]$ identifies
to
\[\prod_{j=0}^i s^j H^{2j-i}(k,A(j))!\;\text{.}\]
\end{proposition3}

For $i\geq 0$, it follows from $A$-valued variants of
corollary~\ref{corollary-morphisms-z-gr-k-z-n-2n-i} and
lemma~\ref{lemma-omega-chi-n}. If $i<0$, $K(A(n),2n-i)$ identifies to
$\R\SheafHom_\bullet(\Gm^{\wedge -i},K(A(n),2n))$ and both
projective systems vanish.

\bigskip

Then, we observe that for any abelian group $A$, $\lim A!\simeq
\Hom(\mathbf{Q},A)$ and $\R^1\lim A!\simeq \Ext(\mathbf{Q},A)$, and that
the shift functor $s$ does not change $\lim$ and $\R^1\lim$ of projective
systems. Thus, we get the following theorem:

\begin{theorem3}
Let $k$ be a perfect field. Let $A$ be an abelian group.
Let $i\in\mathbf{Z}$. There is a canonical short exact sequence:
\begin{multline*}
0\to \prod_{j=0}^{i+1} \Ext(\mathbf{Q},H^{2j-i-1}(k,A(j)))
\to
\Hom_{\SH[k]}(\BGL,\HMot[A][-i])\\
\to \prod_{j=0}^i \Hom(\mathbf{Q},H^{2j-i}(k,A(j)))\to 0\;\text{,}
\end{multline*}
where the group on the right side identifies to morphisms in $\SHnaive[k]$
and the group on the left to stably phantom morphisms.
\end{theorem3}

\begin{corollary3}[Existence of nonzero stably phantom morphisms]
Let $k$ be a perfect field. There exists an isomorphism
\[\Hom_{\SH[k]}(\BGL,\HMot{}[1])\simeq \Ext(\mathbf{Q},\mathbf{Z})\simeq
\widehat{\mathbf{Z}}/\mathbf{Z}\;\text{,}\]
and all these morphisms $f\colon \BGL\to\HMot{}[1]$ are stably phantom,
\ie, for any morphism in $\SH[k]$ of the form $g\colon
(\mathbf{P}^1)^{\wedge -n}\wedge W\to \BGL$ where $n\in\mathbf{Z}$ and
$W\in\Hopt[k]$, the composition $f\circ g$ vanishes
(see \cite[D\'efinition~6.6]{riou-smf}).
\end{corollary3}

\begin{remark3}
Most of the results appearing in this article have homologues in the
classical homotopy theory and are coherent with ``complex points functors''
from $\mathbf{A}^1$-homotopy categories to usual (topological)
homotopy categories.
In particular, the spectrum $\BGL(\mathbf{C})$
obtained as the image of $\BGL$ by the ``complex points functors''
$\SH[\mathbf{C}]\to \SHtop$ represents topological complex $K$-theory
(see \cite[Remarque~2.16]{riou-smf}).
Then, if $\Htop\in\SHtop$ is the classical Eilenberg-Mac Lane spectrum, we
get the same computation of the group
$\Hom_{\SHtop}(\BGL(\mathbf{C}),\Htop{}[1])$. The example of
stably phantom morphisms in $\SHtop$ which we hereby get 
may be considered as simpler
than those constructed by Christensen \cite[Proposition~6.10]{christensen}.
\end{remark3}

\begin{definition3}
Let $k$ be a perfect field. We let $\CH\colon \BGL\to \HMot[\mathbf{Q}]$
be the canonical generator of $\Hom_{\SH[k]}(\BGL,\HMot[\mathbf{Q}])\simeq
\mathbf{Q}$. This
is the Chern character. Using Bott periodicity ($\BGL\simeq
\R\SheafHom_\bullet(\mathbf{P}^1,\BGL)$), we deduce from it a sequence of
morphisms $\CH_i\colon \BGLQ\to \HMot[\mathbf{Q}(i)][2i]$ where
$\HMot[\mathbf{Q}(i)]=\HMot[\mathbf{Q}]\wedge (\mathbf{P}^1)^{\wedge
i}[-2i]$. The total Chern character is $\prod_i \CH_i$:
\[\CH_t\colon \BGLQ\to \prod_{i\in \mathbf{Z}}\HMot[\mathbf{Q}(i)][2i]\]
(the infinite product on the right is also a direct sum).
\end{definition3}

\begin{remark3}Remark~\ref{remark-h-beilinson} may be continued as follows.
One easily sees that the Chern character $\CH\colon \BGLQ\to
\HMot[\mathbf{Q}]$ vanishes on $\BGLQ^{(i)}$ (see 
theorem~\ref{theorem-diagonalisation-bglq}) for $i\neq 0$
so that it factors through its direct factor $\HBeilinson$ (see
definition~\ref{definition-h-beilinson}) as
$\BGLQ\to \HBeilinson\vers{\CH^{(0)}} \HMot[\mathbf{Q}]$. It follows from
known results (see \cite{levine-coniveau}) that $\CH^{(0)}\colon
\HBeilinson\to \HMot[\mathbf{Q}]$ is an isomorphism; equivalently,
$\CH_t\colon \BGLQ\to \oplus_{i\in\mathbf{Z}}\HMot[\mathbf{Q}(i)][2i]$
is an isomorphism.
\end{remark3}

\subsection{Grothendieck-Riemann-Roch}

For simplicity, we only consider the case of a projective and smooth
morphism $f\colon X\to S$ in $\Sm[k]$ where $k$ is a perfect
field. We let
$d$ be the relative dimension of $f$. The ``restriction'' of
$\HMotQ\in\SH[k]$ to $X$ and $S$ provides objects in $\SH[X]$ and $\SH[S]$
which shall also be denoted $\HMotQ$; they satisfy $\L f^\star \HMotQ\simeq
\HMotQ$. To $f$ is attached a morphism of motives $\mathbf{Z}(d)[2d]\to
M(X)$ in $\DM$ (see \cite[\S{}I.4.4]{ivorra-these}) which induces a
morphism $f_\star\colon \R\SheafHom_\bullet(X_+,\HMotQ)\to
\HMot[\mathbf{Q}(-d)][-2d]$ in $\SH$. This morphism induces the pushforward
maps \[f_\star\colon H^p(X\times_S T,\mathbf{Z}(q))\to
H^{p-2d}(T,\mathbf{Z}(q-d))\;\text{,}\] for all $T\in\Sm$.

\begin{theorem}
Let $k$ be a perfect field.
Let $f\colon X\to S$ be a projective and smooth morphism in $\Sm[k]$.
Then, the following diagram commutes in $\SH[S]$:
\[
\xymatrix@C=3cm{\R f_\star \BGLQ\ar[d]^{f_\star}
\ar[r]^-{\R f_\star(\CH\cdot \Td
T_f)} & \prod_{i\in\mathbf{Z}} \R f_\star \HMot[\mathbf{Q}(i)][2i]
\ar[d]^{f_\star}\\
\BGLQ\ar[r]^-{\CH} & \prod_{i\in\mathbf{Z}} \HMot[\mathbf{Q}(i)][2i]}
\]
\end{theorem}

The proof is similar to that of
theorem~\ref{theorem-adams-riemann-roch-sh}. This statement is equivalent
to the usual Grothendieck-Riemann-Roch theorem for morphisms $f_T\colon
X\times_S T\to T$ for all $T\in\Sm[S]$ (which is known to be true, see
\cite[Chapter~15]{fulton}). The reason for this is the variant
of theorem~\ref{theorem-morphisms-f-star-bgl-to-bgl}:
a morphism $\tau\colon \BGLQ\to
\HMot[\mathbf{Q}(i)][2i]$ vanishes if and only if it vanishes after the
application of functors $\pi_{2n,n}\colon \SH \to \Sm^\opp\Ab$
for all $n\in\mathbf{Z}$.

\begin{corollary}
Let $k$ be a perfect field.
Let $f\colon X\to S$ be a projective and smooth morphism in $\Sm[k]$.
For any $j\in\mathbf{N}$, the following diagram commutes:
\[
\xymatrix@C=2cm{
K_j(X)\ar[r]^-{\CH_t\cdot \Td Tf}\ar[d]^{f_\star} & \ar[d]^{f_\star}
\prod_{i\in\mathbf{Z}} H^{2i-j}(S,\mathbf{Q}(i))\\
K_j(S)\ar[r]^-{\CH_t} & \prod_{i\in\mathbf{Z}} H^{2i-j}(S,\mathbf{Q}(i))
}
\]
\end{corollary}

This gives another proof of some results by Gillet \cite{gillet} on higher
Riemann-Roch theorems.



\bibliography{operations}
\bibliographystyle{amsplain}


\end{document}